\documentclass[11pt,reqno,twoside]{amsart}
\usepackage{amssymb,mathrsfs}

\makeatletter
\newif\ifappendix@
\newif\ifex@
\appendix@false
\ex@false

\renewcommand{\c}{\mathfrak c}
\newcommand{\g}{\mathfrak g}
\renewcommand{\r}{\mathfrak r}

\renewcommand{\t}{\mathfrak t}
\renewcommand{\b}{\mathfrak b}
\renewcommand{\l}{\mathfrak l}
\renewcommand{\k}{\mathfrak k}
\newcommand{\z}{\mathfrak z}
\newcommand{\m}{\mathfrak m}
\newcommand{\n}{\mathfrak n}
\renewcommand{\u}{\mathfrak u}
\newcommand{\p}{\mathfrak p}

\newcommand{\s}{\mathfrak s}

\renewcommand{\sl}{\mathfrak{sl}}
\newcommand{\so}{\mathfrak{so}}
\renewcommand{\sp}{\mathfrak{sp}}

\newcommand{\height}{\operatorname{ht}}
\newcommand{\eps}{\varepsilon}
\newcommand{\col}{\operatorname{col}}
\newcommand{\row}{\operatorname{row}}

\newcommand{\C}{\mathbb C}
\newcommand{\R}{\mathbb R}
\newcommand{\Q}{\mathbb Q}
\newcommand{\Z}{\mathbb Z}

\DeclareMathOperator{\ad}{ad}
\DeclareMathOperator{\End}{End}
\DeclareMathOperator{\op}{op}

\newtheorem{thm}{Theorem}
\newtheorem{lem}[thm]{Lemma}
\newtheorem{cor}[thm]{Corollary}

\newtheorem{example}[thm]{Example}

\title[Good grading polytopes]
{\boldmath Good grading polytopes}
\author[J.~Brundan and S.~M.~Goodwin]
{Jonathan Brundan and Simon M.~Goodwin}
\date{}

\address{Department of Mathematics, University of Oregon, Eugene, OR 97403, USA}
\email{brundan@darkwing.uoregon.edu}

\address{Institut for Matematiske Fag, Aarhus Universitet, DK-8000, Aarhus C, Denmark.}
\email{goodwin@imf.au.dk}

\thanks{2000 {\it Mathematics Subject Classification}: 17B20.}
\thanks{First author supported in part by NSF grant no. DMS-0139019.}
\thanks{Second author supported in part by
LIEGRITS and CAALT}

\begin{document}

\begin{abstract}
Let $\g$ be a finite dimensional semisimple Lie algebra over $\C$
and $e \in \g$ a nilpotent element. Elashvili and Kac have recently
classified all {good $\Z$-gradings} for $e$. We instead consider
{\em good $\R$-gradings}, which are naturally parameterized by an
open convex polytope in a Euclidean space arising from the reductive
part of the centralizer of $e$ in $\g$. As an application, we prove
that the isomorphism type of the {\em finite $W$-algebra} attached
to a good $\R$-grading for $e$ is independent of the particular
choice of good grading. \iffalse In this paper we study the {\em
restricted root system} $\Phi_e$ of the centralizer of $e$ in $\g$.
This arises by restricting the root system $\Phi$ of $\g$ to Our
focus is mainly on the geometry of chambers and alcoves in the real
vector space generated by $\Phi_e$. As an application, we explain
how {\em good $\R$-gradings for $e$} are parameterized by an open
convex polytope in this space, generalizing recent work of Elashvili
and Kac which classified good $\Z$-gradings, and then prove that the
isomorphism type of the {\em finite $W$-algebra} attached to a good
$\R$-grading for $e$ is independent of the particular choice of good
grading. \fi
\end{abstract}

\maketitle

\vspace{-5mm}
\section{Introduction}

In this article, we construct isomorphisms between the
finite $W$-algebras
associated to a nilpotent orbit in a complex semisimple Lie algebra.
In some important special cases, these finite $W$-algebras
were first defined and studied in the Ph.D. thesis
of Lynch \cite{Ly}, generalizing
a construction of Kostant \cite{K}.
The same algebras were
later rediscovered by mathematical physicists, who coined the
name ``finite $W$-algebra'' used here; see e.g. \cite{debtji}.
In full generality, a finite $W$-algebra associated to
an arbitrary nilpotent orbit was introduced only recently
by Premet \cite{premet1}, who views the resulting algebra as
an enveloping algebra for
the Slodowy slice through the nilpotent orbit in question;
see also \cite{gangin}.

To review a slight generalization of
Premet's definition in more detail,
let $\g$ be a finite dimensional semisimple Lie algebra over $\C$ and
let $e \in \g$ be nilpotent. An $\R$-grading
$$
\Gamma : \g = \bigoplus_{j \in \R} \g_j
$$
of $\g$ is called a {\em good grading} for $e$ if $e \in \g_2$ and
the linear map $\ad e : \g_j \to \g_{j+2}$ is injective for all $j
\le -1$ and surjective for all $j \ge -1$. This definition
originates in \cite{kacroawak}. We call a good grading {\em
integral} if $\g_j = 0$ for all $j\notin\Z$ and {\em even} if $\g_j
= 0$ for all $j \notin 2 \Z$; these are the most important cases. A
classification of all integral good gradings can be found in
\cite{elakac}. By \cite[Theorem 2.1]{elakac}, even good gradings
correspond to {\em nice parabolic subalgebras} as have been
independently classified by Baur and Wallach \cite{bauwal}.

Suppose $\Gamma$ is a good grading for $e$,
and let $(\,,)$ denote the Killing form on $\g$.
The
alternating bilinear form $\langle\,,\rangle$ on $\g_{-1}$ defined by
$\langle x,y \rangle = ([x,y],e)$
is non-degenerate.  Choose a Lagrangian subspace $\k$ of $\g_{-1}$
and define
$$
\m = \k \oplus \bigoplus_{j < -1} \g_j.
$$
This is a nilpotent subalgebra of $\g$ and the map
$\chi : \m \to \C , x \mapsto (x,e)$
defines a representation of $\m$.
The {\em finite $W$-algebra} associated to $e$ and the good grading
$\Gamma$ may then be defined as the endomorphism algebra
$$
H_\chi = \End_{U(\g)}(U(\g) \otimes_{U(\m)} \C_\chi)^{\op}
$$
of the {\em generalized Gelfand-Graev representation} $U(\g)
\otimes_{U(\m)} \C_\chi$. A critical point for this
article is that the definition of the algebra $H_\chi$ is
independent of the particular choice of the Lagrangian subspace
$\k$. More precisely, given another Lagrangian subspace $\k'$ of
$\mathfrak{g}_{-1}$, a construction due to Gan and Ginzburg
\cite{gangin} yields a canonical isomorphism between the finite
$W$-algebras $H_\chi$ and $H_{\chi'}$ arising from the choices $\k$
and $\k'$, respectively.

As explained in detail in the introduction of \cite{brukle}, the
algebras considered originally by Kostant and Lynch in \cite{K,Ly}
are naturally identified with the algebras $H_\chi$ defined here in
the special case that the good grading $\Gamma$ is even, i.e.\ when
the good grading arises from a nice parabolic subalgebra. In
particular, in the even case, $H_\chi$ can actually be realized as a
subalgebra of $U(\mathfrak{p})$, where $\mathfrak{p}$ is the
parabolic subalgebra $\mathfrak{p}= \bigoplus_{i \geq 0} \g_i$. This
makes the representation theory of $H_\chi$ easier to study in the
even case; for instance, it is clear in these cases that $H_\chi$
possesses many finite dimensional irreducible representations
arising from restrictions of finite dimensional
$U(\mathfrak{p})$-modules. In general it is still an open problem to
show even that $H_\chi$ has a one dimensional representation; see
\cite[Conjecture 3.1]{premet2}.

On the other hand, the algebras studied by Premet
\cite{premet1,premet2} and Gan and Ginzburg \cite{gangin} are the
algebras $H_\chi$ defined here in the special case that the good
grading $\Gamma$ is the {\em Dynkin grading}, i.e.\ the grading
defined by embedding $e$ into an $\mathfrak{sl}_2$-triple $(e,h,f)$
and considering the $\ad h$-eigenspace decomposition of
$\mathfrak{g}$. Representation theory of $\sl_2$ implies that the
Dynkin grading is always an integral good grading for $e$. The
present definition of $H_\chi$, involving an arbitrary choice of
good grading for $e$, gives a general framework containing both the
Kostant-Lynch construction and the Premet construction as special
cases. Our main result shows that in fact the algebra $H_\chi$ only
depends up to isomorphism on $e$, not on the choice of good grading
for $e$.

\begin{thm}\label{thm1}
The finite $W$-algebras
$H_\chi$ and $H_{\chi'}$ associated to
any two good gradings $\Gamma$ and $\Gamma'$
for $e$ are isomorphic.
\end{thm}

To prove the theorem, we need to make precise the physicists' idea
of deforming
one good grading into another. Say two good gradings $\Gamma:\g = \bigoplus_{i \in \R} \g_i$
and $\Gamma':\g = \bigoplus_{j \in \R} \g_j'$
are {\em adjacent} if
$$
\g = \bigoplus_{i^- \leq j \leq i^+} \g_i \cap \g_j',
$$
where $i^-$ denotes the largest integer strictly smaller than $i$ and
$i^+$ denotes the smallest integer strictly greater than $i$.
If $\Gamma$ and $\Gamma'$ are adjacent, then
by Lemma~\ref{lem3} below,
there exist Lagrangian subspaces $\k$ in $\g_{-1}$
and $\k'$ in $\g_{-1}'$ such that
$$
\k \oplus \bigoplus_{i < -1} \g_i = \k' \oplus \bigoplus_{j < -1} \g_j',
$$
i.e.\ the nilpotent subalgebra $\m$ defined from $\Gamma$ and $\k$
coincides with the nilpotent subalgebra $\m'$ defined from $\Gamma'$
and $\k'$. With these choices, the algebras $H_\chi$ and $H_{\chi'}$
corresponding to $\Gamma$ and $\Gamma'$ are simply {\em equal}. In
view of the aforementioned result of Gan and Ginzburg (independence
of choice of Lagrangian subspace), Theorem 1 therefore follows if we
can prove that any two good gradings for $e$ are linked by a chain
of adjacent good gradings. The precise statement is as follows.

\begin{thm}\label{thm2} Given any two good gradings $\Gamma$ and $\Gamma'$ for $e$,
there exists a chain $\Gamma_1, \dots, \Gamma_n$ of good gradings for $e$
such that $\Gamma$ is conjugate to $\Gamma_1$,
$\Gamma_{i}$ is adjacent to $\Gamma_{i+1}$
for each $i=1,\dots,n-1$,
and
$\Gamma_n$ is conjugate to
$\Gamma'$.
\end{thm}

To prove Theorem 2, there is a simple geometric picture. To the
nilpotent element $e$, we will explain how to associate an open
convex polytope $\mathscr P_e$ in a Euclidean space of dimension
equal to the rank of the reductive part of the centralizer of $e$,
together with a finite group $W_e$ of symmetries of $\mathscr P_e$,
in such a way that conjugacy classes of good gradings for $e$ are
parameterized by $W_e$-orbits on $\mathscr P_e$.
We call $\mathscr
P_e$ the {\em good grading polytope}. There is then a natural
collection of affine hyperplanes which cuts the good grading
polytope into finitely many connected {\em alcoves}, with the
property that the good gradings parameterized by points $p$ and $p'$
are adjacent if and only if $p$ and $p'$ lie in the closure of the
same alcove. Since one can get from any point in $\mathscr P_e$ to
any other by crossing finitely many walls, Theorem 2 follows easily
from this description.

The group $W_e$ of symmetries of $\mathscr P_e$
is actually a well known group:
it is isomorphic to the group $N_W(W_J) / W_J$, where
$W_J$ is the parabolic subgroup
of the Weyl group $W$ corresponding
to the minimal Levi subalgebra of $\g$ containing $e$ according to
the Bala-Carter theory. We point out especially
Lemma~\ref{lem11} below which gives
another sense in which these groups
are ``almost'' reflection groups, different to
that of Howlett \cite{How}.
In general, the inequalities defining the good grading polytope
and the hyperplane arrangement defining the alcoves are
all naturally described in terms of what we call the {\em restricted
root system} of the centralizer $\g_e$ of $e$ in $\g$. These
restricted root systems are easy to compute on a case-by-case
basis for the exceptional Lie algebras. For classical Lie algebras,
we instead follow the approach of \cite{elakac} to give a uniform
description of the good grading
polytopes exploiting the natural representation rather than the
adjoint representation.
We expect that the restricted root systems
investigated here will also play a role in the
representation theory of the finite $W$-algebras $H_\chi$
themselves.

\iffalse
In type $A$, it happens that any two integral good gradings for $e$
are linked by a chain of adjacent integral good gradings, so there is really no
motivation to introduce non-integral gradings
in this case. For other types, this is usually not true,
though it always happens that
any two even good gradings are linked by a chain
of adjacent integral good gradings, and any two
integral good gradings are linked by a chain of
adjacent good $\frac{1}{2}\Z$-gradings.
Hence to construct isomorphisms between the finite $W$-algebras
associated to any two integral good gradings for $e$, it is always
enough just to work with $\frac{1}{2}\Z$-gradings.
\fi

\vspace{2mm}

\noindent
{\em Acknowledgements.} We thank Ross Lawther, Gerhard R\"ohrle,
Gary Seitz, Eric Sommers and
Sergey Yuzvinsky
 for help.

\vspace{2mm}

\section{Restricted root systems}\label{s2}

Let $G$ be a semisimple algebraic group over $\C$,
$T$ be a maximal torus and $B$ be a Borel subgroup containing $T$.
We write $\g$, $\t$ and $\b$ for the corresponding Lie algebras.
Recall some standard notation:
\begin{itemize}
\item[-] $\Phi \subset \t^*$ denotes the root system
of $\g$ with respect to $\t$;
\item[-]
$\g_\alpha$ denotes the $\alpha$-root space of $\g$
for each $\alpha \in \Phi$;
\item[-] $\Phi^+ \subseteq \Phi$ is the system of positive roots
defined from $\b = \t \oplus \sum_{\alpha \in \Phi^+} \g_\alpha$;
\item[-] $\Delta =\{\alpha_1,\dots,\alpha_r\}$
is the corresponding set of simple roots;
\item[-] $E$ is the $\R$-lattice
$\R\alpha_1\oplus\cdots\oplus \R \alpha_r$ in $\t^*$
and $E^*$ is the dual lattice
in $\t$;
\item[-] $H_\alpha = \ker \alpha$
is the hyperplane
in $E^*$ defined by $\alpha \in \Phi$;
\item[-]
$\mathscr A$ is the hyperplane arrangement $\{H_\alpha \mid \alpha
\in \Phi\}$ in the real vector space $E^*$;
\item[-] $W < GL(E^*)$ is the Weyl group generated by the simple reflections $s_1,\dots,s_r$, where
$s_i$ is the reflection in the hyperplane $H_i = H_{\alpha_i}$.
\end{itemize}
Given in addition a subset
$J$ of $\{1,\dots,r\}$, we adopt some more standard notation for
parabolic objects associated to $J$:
\begin{itemize}
\item[-] $E_J$ denotes $\sum_{j \in J} \R \alpha_j \subseteq E$;
\item[-] $\Phi_J = \Phi \cap E_J$ is the closed subsystem of
$\Phi$ generated by $\{\pm \alpha_j \mid j \in J\}$ with base
$\Delta_J = \{\alpha_j \mid j \in J\}$;
\item[-] $\p_J = \l_J \oplus \u_J$ denotes the standard parabolic subalgebra
of $\g$ with Levi subalgebra $\l_J = \t \oplus
\sum_{\alpha \in \Phi_J} \g_\alpha$
and nilradical $\u_J = \sum_{\alpha \in \Phi^+ \setminus \Phi_J} \g_\alpha$;
\item[-] $P_J = L_J U_J$ is the corresponding standard parabolic subgroup of $G$ with standard Levi subgroup $L_J$ and unipotent radical $U_J$;
\item[-] $W_J$ denotes the parabolic subgroup of $W$ generated
by $\{s_j \mid j \in J\}$.
\end{itemize}
Our final piece of notation is less standard:
$E^J \cong (E / E_J)^*$ denotes
$\bigcap_{j \in J} H_j \subseteq E^*$.
Then,
$$
\mathscr A^J = \{H_\alpha \cap E^J \mid \alpha \in \Phi \setminus
\Phi_J\}
$$
is the restriction of the reflection arrangement $\mathscr A$ to the
subspace $E^J$. It has been well studied in the literature, starting
from work of Orlik and Solomon \cite{OS}. For $\alpha \in E$, we let
$\alpha^J \in E/ E_J$ denote the restriction of $\alpha$ to $E^J$.
In this section, we want to focus not on the restricted arrangement
$\mathscr A^J$, but rather on the {\em restricted root system}
$$
\Phi^J = \{\alpha^J  \mid \alpha \in \Phi \setminus \Phi_J\}
$$
consisting of all the non-zero restrictions of roots in $\Phi$ to
$E^J$.  The hyperplanes in $\mathscr A^J$ are the kernels of the
restricted roots in $\Phi^J$, so one can recover $\mathscr A^J$ from
$\Phi^J$, but not vice versa. Note that $\Phi^J$ is in general
definitely {\em not} a root system in $E / E_J$ in the usual sense.

%We note here that there is a related construction of the {\em
%quotient root system} of $\Phi$ by $\Phi_J$, see \cite[10.4]{C}. The
%quotient root system is a subset of $\Phi_J$ that is a root system
%in the subspace of $E / E_J$ that it spans.

%It seems likely that the quotient root system coincides with the
%root system in $\Phi_e^\circ$, from Section \ref{s3} of this paper, in case
%$e$ is regular nilpotent in $L_J$.

From now on, we will always identify $E$ with $E^*$ using the real
inner product $(.,.)$ induced by the Killing form on $\g$.
We can then identify both the spaces $E^J$ and $E / E_J$
with the orthogonal complement to $E_J$ in $E$.
Under this identification, the notation $\alpha^J$
becomes the orthogonal projection of $\alpha \in E$ to $E^J$
along the direct sum decomposition $E = E^J \oplus E_J$.
Let us also set
$I = \{1,\dots,r\} \setminus J$ and $m = |I| = \dim E^J$.

\begin{lem}
For any $\alpha \in \Phi^J$, there exists $\alpha' \in E_J$
such that $\alpha+\alpha' \in \Phi$
and $(\alpha', \alpha_j) \geq 0$ for all $j \in J$.
\end{lem}

\begin{proof}
By the definition of $\Phi^J$, the set $\{\alpha' \in E_J \mid
\alpha+\alpha' \in \Phi\}$ is non-empty. Pick an element $\alpha'$
from this set that is maximal in the dominance ordering. To complete
the proof, we just need to show that $(\alpha',\alpha_j) \geq 0$ for
all $j \in J$. Well, if not, we can find $j \in J$ such that
$(\alpha+\alpha',\alpha_j) = (\alpha',\alpha_j) < 0$, but then
$\alpha + (\alpha'+\alpha_j) \in \Phi$ by \cite[Lemma 9.4]{Hum}
contradicting the maximality of the choice of $\alpha'$.
\end{proof}

\begin{lem}\label{prop1}
If $\alpha,\beta \in \Phi^J$ are distinct
roots with $(\alpha,\beta) > 0$, then $\alpha-\beta \in \Phi^J$ too.
\end{lem}

\begin{proof}
By the previous lemma, we can lift $\alpha,\beta$ to
$\alpha+\alpha', \beta+\beta' \in \Phi$, where $\alpha',\beta' \in E_J$
satisfy $(\alpha',\alpha_j) \geq 0, (\beta',\alpha_j) \geq 0$ for all $j \in J$.
In other words, $\alpha', \beta'$ belong to the closure of the same chamber
in $E_J$, hence $(\alpha',\beta') \geq 0$.
So
$$
(\alpha+\alpha',\beta+\beta') = (\alpha,\beta)+(\alpha',\beta') > 0.
$$
Also $\alpha+\alpha' \neq \beta+\beta'$ since $\alpha \neq \beta$.
So \cite[Lemma 9.4]{Hum} implies that
$(\alpha+\alpha')-(\beta+\beta') \in \Phi$.
Hence, $\alpha-\beta \in \Phi^J$.
\end{proof}

\begin{lem}\label{prop2}
If $\alpha,\beta \in\Phi^J$ are proportional roots,
then there exists $\gamma \in\Phi^J$ such that $\alpha$, $\beta$
are both integer multiples of $\gamma$.
\end{lem}

\begin{proof}
Let $M = \{c > 0  \mid c \alpha \in \Phi^J\}$. The previous lemma
implies that if $a,b$ are distinct elements of $M$ then $|a-b| \in
M$ too. It follows that any element of $M$ is an integer multiple of
the smallest element.
\end{proof}

Define a {\em base} of the restricted root system $\Phi^J$ to be a
subset $\{\beta_i \mid i \in I\}$ of $\Phi^J$ such that any element
of $\Phi^J$ can be written as $\sum_{i \in I} a_i \beta_i$ with
either all $a_i \in \Z_{\geq 0}$ or all $a_i \in \Z_{\leq 0}$. Of
course any base for $\Phi^J$ is necessarily a basis for the vector
space $E^J$. Any base $\{\beta_i \mid i \in I\}$ partitions the
restricted root system $\Phi^J$ into positive and negative roots,
the positive ones being the roots that are a positive linear
combination of $\beta_i$'s. To construct bases of $\Phi^J$ in the
usual way, let $\gamma \in E^J$ be {\em regular}. This means that
$\gamma$ does not lie on any of the hyperplanes in $\mathscr A^J$,
or equivalently, $(\alpha, \gamma) \neq 0$ for all $\alpha \in
\Phi^J$. Then we can define $\Phi^J(\gamma)$ to be $\{\alpha \in
\Phi^J \mid  (\alpha, \gamma) > 0\}$, and clearly $\Phi^J =
\Phi^J(\gamma) \sqcup (- \Phi^J(\gamma))$. Call $\alpha \in
\Phi^J(\gamma)$ {\em decomposable} if $\alpha = \beta_1+\beta_2$ for
$\beta_1,\beta_2 \in \Phi^J(\gamma)$, and {\em indecomposable}
otherwise. Armed with Lemmas~\ref{prop1} and \ref{prop2}, the
following theorem is proved in essentially the same way as for root
systems, see e.g. \cite[Theorem 10.1]{Hum}.

\begin{thm}\label{bases}
Let $\gamma \in E^J$ be regular. Then the set
$\Delta^J(\gamma)$ of all indecomposable roots in
$\Phi^J(\gamma)$
is a base of
$\Phi^J$, and every base can be obtained in this manner.
\end{thm}

Recall that the bases for the root system $\Phi$ are in natural
bijective correspondence with set $\mathscr C$ of {\em chambers} in
the hyperplane arrangement $\mathscr A$, that is, the connected
components of $E \setminus \bigcup \mathscr A$. Under this
correspondence, the base $\{\beta_1,\dots,\beta_r\}$ corresponds to
the chamber $\{\alpha \in E \mid (\alpha,\beta_i) > 0\text{ for all
}i=1,\dots,r\}$. Theorem~\ref{bases} leads to a similar bijection
between the set of bases of the restricted root system $\Phi^J$ and
the set $\mathscr C^J$ of chambers in the hyperplane arrangement
$\mathscr A^J$.

\begin{cor}\label{oto}
There is a natural bijective correspondence between bases in
$\Phi^J$ and chambers in $\mathscr C^J$, under which the base
$\{\beta_i \mid i \in I\}$ corresponds to the chamber $\{\alpha \in
E^J \mid (\alpha,\beta_i) > 0 \text{ for all }i \in I\}$.
\end{cor}

\begin{proof}
We just explain how to construct the inverse map from chambers to bases.
Given a chamber $C \in \mathscr C^J$, pick any (necessarily regular) point
$\gamma \in C$. Then, image of $C$ under the inverse map is
the base $\Delta^J(\gamma)$ for $\Phi^J$.
This is well-defined, because if $\gamma, \gamma'$
belong to the same chamber, then
they
lie on the same side of each hyperplane in $\mathscr A^J$, so
$\Phi^J(\gamma) = \Phi^J(\gamma')$.
\end{proof}

There is another way to construct bases for the restricted root system
$\Phi^J$, by restricting
bases for $\Phi$ that contain bases for $\Phi_J$.

\begin{lem}\label{p1}
Suppose that $\{\beta_1,\dots,\beta_r\}$ is a base for $\Phi$ such
that $\{\beta_j \mid j \in J\}$ is a base for $\Phi_J$. Then,
$\{\beta_i^J \mid i \in I\}$ is a base for $\Phi^J$, and every base
for $\Phi^J$ can be obtained in this way.
\end{lem}

\begin{proof}
Suppose first that $\{\beta_1,\dots,\beta_r\}$ is a base for $\Phi$
such that $\{\beta_j \mid j \in J\}$ is a base for $\Phi_J$. Any
$\alpha \in \Phi \setminus \Phi_J$ can be written as $\alpha =
\sum_{j=1}^r a_j \beta_j$, so that the $a_j$'s are either all $\geq
0$ or all $\leq 0$. Since $\beta_j \in E_J$ for each $j \in J$,
$\alpha^J = \sum_{i \in I} a_i \beta_i^J$. Hence, $\{\beta_i^J \mid
i \in I\}$ is a base for $\Phi^J$.

To show that every base in $\Phi^J$ arises in this way, we think
instead in terms of chambers. Let $C \in \mathscr C$ be the chamber
corresponding to the base $\{\beta_1,\dots,\beta_r\}$, still
assuming that $\{\beta_j \mid j \in J\}$ is a base for $\Phi_J$. The
closure $\overline{C}$ is equal to $\{\alpha \in E \mid
(\alpha,\beta_i) \geq 0 \text{ for all }i=1,\dots,r\}$, while $E^J =
\{\alpha \in E \mid (\alpha,\beta_j) = 0 \text{ for all }j \in J\}$.
Hence, the intersection $\overline{C} \cap E^J$ is equal to
$\{\alpha \in E^J \mid (\alpha,\beta_i^J) \geq 0 \text{ for all }j
\in J\}$. This shows that $(\overline{C} \cap E^J) \setminus \bigcup
\mathscr A^J$ is the  chamber in $\mathscr C^J$ corresponding to the
base $\{\beta_i^J \mid i \in I\}$. We must prove that every chamber
in $\mathscr C^J$ can be obtained in this way.

Suppose that $\{\beta_1,\dots,\beta_r\}$ is a base
for $\Phi$ that does {\em not} contain a base for $\Phi_J$, and let
$C$ be the corresponding chamber in $\mathscr C$.
We can find $\beta = \sum_{j=1}^r a_j \beta_j \in \Phi_J$ such that
$a_i \neq 0$ for some $1 \leq i \leq r$ with $\beta_i \notin \Phi_J$.
Take any $\alpha \in \overline{C} \cap E^J$, so
$(\alpha,\beta) = 0$ and
$(\alpha,\beta_j) \geq 0$ for all $j=1,\dots,r$.
Since $a_i \neq 0$, the equation $\sum_{j=1}^r a_j (\alpha,\beta_j) = 0$
implies that $(\alpha,\beta_i) = 0$.
Hence, $\overline{C} \cap E^J$ is contained in the hyperplane
$H_{\beta_i}$, and
$(\overline{C} \cap E^J) \setminus \bigcup \mathscr A^J = \varnothing.$
Since $E^J \setminus \bigcup \mathscr A^J$ is obviously covered by the
sets $(\overline C \cap E^J) \setminus \bigcup  \mathscr A^J$
as $C$ runs over all chambers in $\mathscr A$,
we have now shown that every chamber in
$\mathscr C^J$ is equal
to $(\overline{C} \cap E^J) \setminus \bigcup  \mathscr A^J$
for some chamber $C$ in $\mathscr C$ such that the corresponding
base of $\Phi$ contains a base for $\Phi_J$.
\end{proof}

\begin{lem}\label{p2}
Suppose that $\{\beta_1,\dots,\beta_r\}$ and
$\{\gamma_1,\dots,\gamma_r\}$ are two bases for $\Phi$ such that
$\{\beta_j \mid j \in J\}$ and $\{\gamma_j \mid j \in J\}$ are bases
for $\Phi_J$. The resulting bases $\{\beta_i^J \mid i \in I\}$ and
$\{\gamma_i^J \mid i \in I\}$ for $\Phi^J$ are equal if and only if
there exists $w \in W_J$ mapping $\{\beta_1,\dots,\beta_r\}$ to
$\{\gamma_1,\dots,\gamma_r\}$.
\end{lem}

\begin{proof}
Since $W_J$ acts trivially on $E^J$, it is easy to see that if
$\{\beta_1,\dots,\beta_r\}$ and $\{\gamma_1,\dots,\gamma_r\}$ are
conjugate under $W_J$, then $\{\beta_i^J \mid i \in I\}$ and
$\{\gamma_i^J \mid i \in I\}$ are equal. Conversely, suppose that
$\{\beta_i^J \mid i \in I\}$ and $\{\gamma_i^J \mid i \in I\}$ are
equal. Recalling that $W_J$ acts transitively on bases for $\Phi_J$,
we can conjugate and reindex if necessary to assume that $\beta_j =
\gamma_j$ for all $j \in J$ and that $\beta_i^J = \gamma_i^J$ for
all $i \in I$. But then we can certainly write
$$
\beta_i = \gamma_i + \sum_{j \in J} a_{i,j} \gamma_j
$$
for every $i \in I$ and scalars $a_{i,j} \in \R$.
Since $\beta_i$ is a root and the $\gamma_i$'s form a base
for $\Phi$, we get that $a_{i,j} \geq 0$ for all $i \in I, j \in J$.
But also
$$
\gamma_i = \beta_i - \sum_{j \in J} a_{i,j} \beta_j
$$
for every $i \in I$,
which implies that all $a_{i,j} \leq 0$ too.
Hence, $\beta_i = \gamma_i$ for each $i \in I$,
and the original bases for $\Phi$ are equal as required.
\end{proof}

\begin{thm}\label{bij}
There is a natural bijective correspondence between bases for $\Phi$
containing $\Delta_J$ and bases for $\Phi^J$, under which the base
$\{\beta_i,\alpha_j \mid i \in I, j\in J\}$ for $\Phi$ corresponds
to the base $\{\beta_i^J \mid i \in I\}$ for $\Phi^J$.
\end{thm}

\begin{proof}
Since $W_J$ acts simply transitively on bases for $\Phi_J$,
each $W_J$-orbit of bases $\{\beta_1,\dots,\beta_r\}$
for $\Phi$ containing a base for $\Phi_J$
has a unique representative that contains $\Delta_J$.
Given this, the theorem is immediate from Lemmas~\ref{p1} and \ref{p2}.
\end{proof}

We remark that bases for $\Phi$ containing $\Delta_J$ are also in
bijective correspondence with parabolic subgroups $P$ of $G$ that
have $L_J$ as a Levi factor, the base $\{\beta_i, \alpha_j \mid i
\in I, j \in J\}$ for $\Phi$ corresponding to the parabolic subgroup
with Lie algebra generated by $\l_J$ and all $\g_{\beta_i}$ ($i \in
I$). So another way of thinking about Theorem~\ref{bij} is that
choosing a base for the restricted root system $\Phi^J$ is
equivalent to choosing a parabolic subgroup $P$ of $G$ with Levi
factor $L_J$, just as choosing a base for $\Phi$ is equivalent to
choosing a Borel subgroup of $G$ containing $T$.

Corresponding to the base $\Delta$ of $\Phi$,
or to the standard parabolic subgroup $P_J$ of $G$,
we have the {\em standard base}
$$
\Delta^J = \{\alpha^J \mid \alpha \in \Delta \setminus \Delta_J\} =
\{\alpha_i^J \mid i \in I\}
$$
of $\Phi^J$.
Now suppose that $K$ is a subset of $\{1,\dots,r\}$ such that
$w \cdot \Delta_K = \Delta_J$
for some $w \in W$.
Since $w \cdot E_K = E_J$, $w$ induces an isometry
between $E^K$ and $E^J$ which maps $\Phi^K$ to $\Phi^J$.
So if we apply $w$ to the standard base $\Delta^K$ of $\Phi^K$,
we obtain a base
$w \cdot \Delta^K$ for $\Phi^J$.
Clearly, all bases for $\Phi$ containing $\Delta_J$ are of the form
$w \cdot \Delta$ for some $K \subset \{1,\dots,r\}$ and some $w \in W$
such that $w \cdot \Delta_K = \Delta_J$.
Therefore, by Theorem~\ref{bij}, all bases for $\Phi^J$ are of the form
$w \cdot \Delta^K$ for suitable $w$ and $K$. This means that for most
purposes, it is sufficient to work only with standard
bases $\Delta^J$, providing one is prepared to
allow the subset $J$ of $\{1,\dots,r\}$ to change.

Finally, we introduce the {\em restricted Weyl group} $W^J$, namely,
the pointwise stabilizer in $W$ of the set $\Delta_J$. This is a
well known group, studied in particular by Howlett \cite{How};
see also \cite[$\S$10.4]{C}.
Clearly, $W^J$ normalizes $W_J$ and $W^J \cap W_J = \{1\}$. In fact,
by \cite[Lemma 2]{How}, we have that $W_J W^J = N_W(W_J)$, so $W^J
\cong N_W(W_J) / W_J$. By Lemma~\ref{enum} below, the natural action
of $W^J$ on $E^J$ is faithful, so we can view $W^J$ as a subgroup of
$GL(E^J)$. In general, $W^J$ is {\em not} a reflection group, though
it is close to being one in a sense made precise in Howlett's work;
we will give an alternative explanation of this phenomenon in the
next section. Clearly, $W^J$ leaves the subset $\Phi^J \subset E^J$
invariant, hence we get an induced action of $W^J$ on the set of
bases for the root system $\Phi^J$.

For the next lemma we require the following piece of notation:
define $\mathscr{K}_J$ to be the set of subsets $K$ of $\{1,\dots,r\}$
with the property that $w \cdot \Delta_K = \Delta_J$ for some $w \in W$.

\begin{lem}\label{enum}
For each $K \in \mathscr{K}_J$, pick $w_K \in W$ such that $w_K\cdot
\Delta_{K} = \Delta_J$. Then,
$$
\{w_K\cdot \Delta^K \mid K \in \mathscr{K}_J\}
$$
is a set of orbit representatives for the action of
the restricted Weyl group $W^J$ on the set of bases for $\Phi^J$.
Moreover, each orbit is regular, of size $|W^J|$.
\end{lem}

\begin{proof}
The set of all $w \in W$ with the property that $\Delta_J \subseteq
w\cdot\Delta$ is the disjoint union $\bigcup_{K \in \mathscr{K}_J}
W^J w_K$. Since $W$ acts simply transitively on bases for $\Phi$,
this means that there are $|\mathscr{K}_J||W^J|$ different bases for $\Phi$
containing $\Delta_J$, namely, the bases $\{w w_K\cdot\Delta \mid w
\in W^J, K \in \mathscr{K}_J\}$. Applying Theorem~\ref{bij}, we
deduce that there are $|\mathscr{K}_J||W^J|$ different bases for $\Phi^J$,
namely, the bases $\{w w_K \cdot \Delta^K \mid w \in W^J, K \in
\mathscr{K}_J\}$. The lemma follows.
\end{proof}

This lemma immediately implies that the number of bases for the
restricted root system $\Phi^J$ is equal to $|\mathscr K_J||W^J|$.
Equivalently, by Corollary~\ref{oto},
the number of chambers in the hyperplane arrangement $\mathscr A^J$
is given by the formula
$$
\left|\mathscr C^J\right| = |\mathscr K_J| \left|W^J\right|.
$$
This is a well known identity due originally to Orlik and Solomon
\cite[(4.2)]{OS}. The hyperplane arrangement $\mathscr A^J$ is known
to be a free arrangement; see \cite{OT,Doug}. So by
\cite[$\S$4.6]{OT2} its Poincar\'e polynomial can be expressed as
$(1+b_1^J t)\cdots (1+b_m^J t)$ for {\em exponents} $b_1^J \leq
\dots \leq b_m^J$. This factorization already appears in \cite{OS},
and the exponents were computed there too in all cases. It is well
known that
\begin{align*}
\left|\mathscr A^J\right| &= b_1^J+b_2^J+\cdots+b_m^J,\\
\left|\mathscr C^J\right| &= (1+b_1^J) (1+b_2^J)\cdots (1+b_m^J).
\end{align*}
Moreover,  if $G$ is simple and $m \geq 1$
then
the arrangement $\mathscr A^J$ is irreducible,
hence $b_1^J = 1$ and,
assuming $m \geq 2$ too, $b_2^J \geq 2$ .

To conclude the section, we want to mention a theorem of Sommers \cite{S}
which gives a quick way to determine the exponents
$b_i^J$.
For any $\alpha = \sum_{i=1}^r a_i \alpha_i \in E$, we let
$\height(\alpha)$ denote $\sum_{i=1}^r a_i$.
Let $\theta = \sum_{i=1}^r c_i \alpha_i$
be the highest root in $\Phi$, and recall that
all other roots in $\Phi$ are strictly smaller than $\theta$ in the dominance
ordering. It follows easily that
$\theta^J$ is the unique highest root in $\Phi^J$.
Now
introduce the following plausible analogue of the Coxeter number for the restricted
root system $\Phi^J$: let
$$
h^J = \min \{\height(\theta^K)+1 \mid K \in \mathscr{K}_J \}.
$$
Then, Sommers' theorem says that an integer $1 \leq p < h^J$ belongs
to the set $\{b_1^J, \dots, b_m^J\}$ of exponents whenever it is
prime to all the coefficients $c_1,\dots,c_r$ of $\theta$. Combined
with the facts mentioned in the previous paragraph, and also
\cite{How} (or Lemmas~\ref{lem1}--\ref{lem11} below) from which the
orders of the groups $W^J$ can be computed, this always give enough
information to completely determine the exponents.

\ifappendix@
For convenience, we have listed the basic numerical invariants of $\Phi^J$
for each exceptional group in Tables \ref{tab1}--\ref{tab2}, using
the notation of \cite{C} for conjugacy classes of Levi subgroup.
Restricted root systems arising from
classical groups will be discussed in more detail later in the article;
see also \cite{OS, S}.
For an example of a rank 2 restricted root system, see Example~\ref{eg1}.
\fi

\begin{example}\label{eg1}\rm
Take $G = E_7$.
Label the simple roots/the vertices of the
Dynkin diagram as follows:
$$
\begin{array}{lllllll}
3&4&2&5&6&7\\
&&1
\end{array}
$$
Take $I = \{1,2\}$ and $J = \{3,4,5,6,7\}$, so $L_J$
is of type $A_3+A_2$.
The positive roots in $\Phi^J$ corresponding to the standard
base $\Delta^J = \{\alpha_1^J,\alpha_2^J\}$
are
$$
\{\alpha_1^J, \alpha_2^J, \alpha_1^J+\alpha_2^J,
\alpha_1^J+2\alpha_2^J, \alpha_1^J+3\alpha_2^J,
2\alpha_1^J+3\alpha_2^J, 2\alpha_1^J+4\alpha_2^J\}.
$$
The restricted Cartan matrix with $ij$-entry
$\frac{2(\alpha_i^J, \alpha_j^J)}{(\alpha_j^J, \alpha_j^J)}$
for $i,j=1,2$
is the matrix $\left(\begin{array}{rr}2 & -\frac{24}{7}\\ -1&2\;\end{array}\right)$.
It follows that $\alpha_1^J \perp (\alpha_1^J+2\alpha_2^J)$.
We get the following picture of roots and orthogonal hyperplanes:
$$
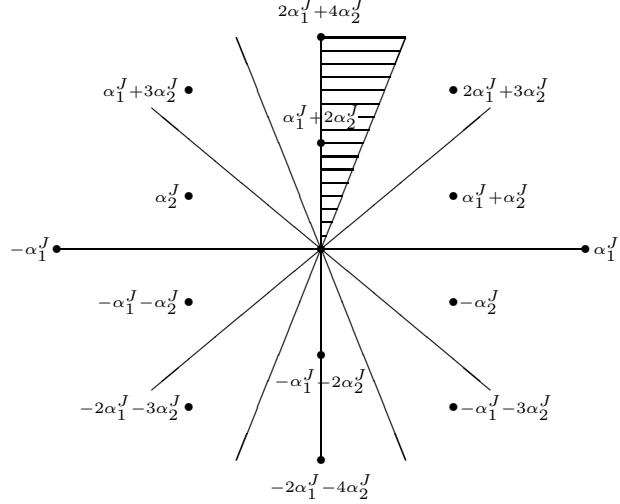
\begin{picture}(200,170)
%\put(0,80){\line(1,0){200}}
%\put(100,00){\line(0,1){160}}
%\put(51,60){\line(5,2){100}}
%\put(51,100){\line(5,-2){100}}
%\put(49.3333,19.2){\line(5,6){100.2}}
%\put(150.6666,19.2){\line(-5,6){100.2}}

\put(100,80){\line(1,0){100}}
\put(100,80){\line(-1,0){100}}
\put(100,80){\line(0,1){80}}
\put(100,80){\line(0,-1){80}}
\put(100,80){\line(-2,5){32}}
\put(100,80){\line(2,-5){32}}
\put(100,80){\line(2,5){32}}
\put(100,80){\line(-2,-5){32}}
\put(100,80){\line(-6,5){64}}
\put(100,80){\line(6,5){64}}
\put(100,80){\line(6,-5){64}}
\put(100,80){\line(-6,-5){64}}

\put(100,85){\line(1,0){2}}
\put(100,90){\line(1,0){4}}
\put(100,95){\line(1,0){6}}
\put(100,100){\line(1,0){8}}
\put(100,105){\line(1,0){10}}
\put(100,110){\line(1,0 ){12}}
\put(100,115){\line(1,0 ){14}}
\put(100,120){\line(1,0){16}}
\put(100,125){\line(1,0){18}}
\put(114,130){\line(1,0){6}}
\put(100,135){\line(1,0){22}}
\put(100,140){\line(1,0){24}}
\put(100,145){\line(1,0){26}}
\put(100,150){\line(1,0){28}}
\put(100,155){\line(1,0){30}}
\put(100,160){\line(1,0){32}}

\put(200,80){\circle*{3}}
\put(0,80){\circle*{3}}
\put(100,0){\circle*{3}}
\put(100,40){\circle*{3}}
\put(100,80){\circle*{3}}
\put(100,120){\circle*{3}}
\put(100,160){\circle*{3}}
\put(150,100){\circle*{3}}
\put(50,60){\circle*{3}}
\put(150,140){\circle*{3}}
\put(50,140){\circle*{3}}
\put(150,60){\circle*{3}}
\put(50,20){\circle*{3}}
\put(150,20){\circle*{3}}
\put(50,100){\circle*{3}}
\put(100,-10){\makebox(0,0){{$\scriptscriptstyle -2\alpha_1^J-4\alpha_2^J$}}}
\put(100,30){\makebox(0,0){{$\scriptscriptstyle -\alpha_1^J-2\alpha_2^J$}}}
\put(100,170){\makebox(0,0){{$\scriptscriptstyle 2\alpha_1^J+4\alpha_2^J$}}}
\put(100,130){\makebox(0,0){{$\scriptscriptstyle \alpha_1^J+2\alpha_2^J$}}}
\put(208,80){\makebox(0,0){{$\scriptscriptstyle \alpha_1^J$}}}
\put(-10,80){\makebox(0,0){{$\scriptscriptstyle -\alpha_1^J$}}}
\put(31,60){\makebox(0,0){{$\scriptscriptstyle -\alpha_1^J-\alpha_2^J$}}}
\put(160,60){\makebox(0,0){{$\scriptscriptstyle -\alpha_2^J$}}}
\put(42,100){\makebox(0,0){{$\scriptscriptstyle \alpha_2^J$}}}
\put(166,100){\makebox(0,0){{$\scriptscriptstyle \alpha_1^J+\alpha_2^J$}}}
\put(32,140){\makebox(0,0){{$\scriptscriptstyle \alpha_1^J+3\alpha_2^J$}}}
\put(170,140){\makebox(0,0){{$\scriptscriptstyle 2\alpha_1^J+3\alpha_2^J$}}}
\put(28,20){\makebox(0,0){{$\scriptscriptstyle -2\alpha_1^J-3\alpha_2^J$}}}
\put(170,20){\makebox(0,0){{$\scriptscriptstyle -\alpha_1^J-3\alpha_2^J$}}}

\end{picture}
$$
\vspace{1mm}

\noindent
There are 12 chambers in the hyperplane arrangement $\mathscr A^J$,
the one corresponding to the standard base
 being shaded.
Since $|\mathscr K_J| = 3$, there are three $W^J$-orbits on
chambers. In fact, $W^J\cong S_2 \times S_2$ is generated by the
reflections in the horizontal and vertical axes. The highest root
$\theta^J$ is $2\alpha_1^J+4\alpha_2^J$, but the Coxeter number
$h^J$ is $6$ not $7$: it comes from the highest root
$2\alpha_1^J+3\alpha_2^J$ with respect to the non-standard base
$\{\alpha_1^J+3\alpha_2^J, -\alpha_2^J\}$. The exponents are $1$ and
$5$.
\end{example}

\section{Centralizers}\label{s3}

We fix for the remainder of the article a nilpotent element $e \in
\g$; our basic references for all matters concerning nilpotent
orbits are \cite[ch.\ 1--5]{Jantzen} and \cite[ch.\ 5]{C}. We denote
the centralizer of $e$ in $G$ either by $Z_G(e)$ or by $G_e$ for
short. Similarly, we write $\z_\g(e)$ or $\g_e$ for its centralizer
in $\g$. By the Jacobson-Morozov theorem, we can embed $e$ into an
$\sl_2$-subalgebra $\s = \langle e,h,f \rangle$, so that $[h,e] =
2e, [h,f] = -2f$ and $[e,f] =h$. Moreover, by a result of Kostant,
any other such triple $(e,h',f')$ is conjugate to $(e,h,f)$ by an
element of the connected centralizer $G_e^\circ$.

The $\ad h$-eigenspace decomposition of $\g$ defines a $\Z$-grading
$\g = \bigoplus_{j \in \Z} \g_j$ which we call the {\em Dynkin
grading}. Let $\c = \g_0$ and let $C$ be the corresponding closed
connected subgroup of $G$. In other words, $\c$ and $C$ are the
centralizers of $h$ in $\g$ and $G$, respectively. Also let $\r =
\bigoplus_{j > 0} \g_j$ and let $R$ be the corresponding closed connected
subgroup of $G$.
It is well known that $C_e$ is a maximal reductive subgroup of
$G_e$, with Lie algebra $\c_e$, and that $R_e$ is the unipotent
radical of $G_e$, with Lie algebra $\r_e$. Moreover, $G_e$ is the
semidirect product $C_e \ltimes R_e$, and $\g_e$ is the semidirect
sum $\c_e \oplus \r_e$. Finally, the component group $G_e /
G_e^\circ$ is isomorphic to $C_e / C_e^\circ$.

Fix a maximal torus $T$ of $G$ contained in $C$ and containing a
maximal torus of $C_e$. An important role is played by the
centralizer $\t_e$ of $e$ in the Lie
algebra $\t$ of $T$. It is a Cartan subalgebra
of the reductive part $\c_e$ of the centralizer $\g_e$. Let $L$ be
the centralizer of $\t_e$ in $G$, and let $\l$ be the Lie algebra of
$L$, i.e.\ the centralizer of $\t_e$ in $\g$. Thus, $L$ is Levi
subgroup of $G$, and the center of $\l$ is equal to $\t_e$. By the
Bala-Carter theory, $\l$ is a minimal Levi subalgebra of $\g$
containing $e$, and $e$ is a distinguished nilpotent element
of the derived subalgebra $[\l,\l]$
of $\l$. Moreover, both  $h$ and $f$ automatically lie in $[\l,\l]$.

\begin{lem}\label{wts} The set
of weights of $\t_e$ on $\g_e$ is equal to
the set of weights of $\t_e$ on $\g$.
\end{lem}

\begin{proof}
For $\alpha \in \t_e^*$ and $i \geq 0$, let
$L(\alpha,i)$ denote the
irreducible $\t_e \oplus \s$-module
of dimension $(i+1)$ on which $\t_e$ acts by weight $\alpha$.
Decompose $\g$ as a $\t_e\oplus \s $-module
$$
\g \cong \bigoplus_{\alpha \in t_e^*} \bigoplus_{i \geq 0}
m(\alpha,i) L(\alpha,i)
$$
for multiplicities $m(\alpha,i) \geq 0$. The set of weights of
$\t_e$ on $\g$ is $\{\alpha\in \t_e^* \mid m(\alpha,i) \neq 0 \text{
for some }i \geq 0\}$. Since each $L(\alpha,i)$ contains a non-zero
vector annihilated by $e$, this is also the set of weights of $\t_e$
on $\g_e$.
\end{proof}

We define $\Phi_e \subset \t_e^*$ to be the set of all non-zero
weights of $\t_e$ on $\g_e$. The zero weight space of $\t_e$ on
$\g_e$ is of course the centralizer $\l_e$ of $e$ in the Levi
subalgebra $\l$. So we have the following analogue of the Cartan
decomposition for centralizers:
$$
\g_e = \l_e \oplus \bigoplus_{\substack{\alpha \in \Phi_e \\ i \geq 0}}
\g_e(\alpha,i)
$$
where $\g_e(\alpha,i) = \{x \in \g_e \mid [h,x] = ix\text{ and
}[t,x] = \alpha(t) x\text{ for all }t \in \t_e\}$. This
decomposition is compatible with the decomposition $\g_e = \c_e
\oplus \r_e$; indeed, we have that
$$\c_e
= \t_e \oplus \bigoplus_{\alpha \in \Phi_e^\circ}
\g_e(\alpha,0),\qquad \r_e  =[\l,\l]_e \oplus
\bigoplus_{\substack{\alpha \in \Phi_e\\ i
> 0}} \g_e(\alpha,i),
$$
where $\Phi_e^\circ$ denotes the set of all $\alpha \in \Phi_e$ such
that $\g_e(\alpha,0)$ is non-zero.  The root system $\Phi_e$ is a
restricted root system in the sense of the previous section. To
explain this, we need to make one more important choice: let $P$ be
a parabolic subgroup of $G$ such that $L$ is a Levi factor of $P$.
Let $U$ be the unipotent radical of $P$, and denote the
corresponding Lie algebras by $\p$ and $\u$, so $\p = \l \oplus \u$.
Pick a Borel subgroup $B$ of $G$ contained in $P$ and containing
$T$. The choices of $T$ and $B$ determine a root system $\Phi$ for
$\g$ and a base $\Delta = \{\alpha_1,\dots,\alpha_r\}$, and we can
appeal to the setup from the previous section. We then have that $L
= L_J$, $U = U_J$ and $P = P_J$ for a unique subset $J$ of
$\{1,\dots,r\}$. The Euclidean space $E^J$ from the previous
section, henceforth denoted $E_e$, is an $\R$-form for the center
$\t_e$ of the Lie algebra $\l = \l_J$, and Lemma~\ref{wts} shows
that $\Phi_e$ coincides with the restricted root system $\Phi^J
\subseteq \t_e^*$. The standard base $\Delta^J$ for $\Phi^J$ will be
denoted from now on by $\Delta_e$, and we let $\Phi_e^+$ denote the
corresponding set of positive roots. In the above root space
decomposition of $\g_e$, we have that $\p_e = \l_e \oplus \u_e$
where
$$
\u_e = \bigoplus_{\substack{\alpha \in \Phi_e^+\\i \ge 0}}
\g_e(\alpha,i).
$$
Indeed, as explained in the previous section, the choice of
the parabolic $P = LU$ is actually
equivalent to this choice $\Phi_e^+$ of positive roots in $\Phi_e$.

The objects $\Phi_e, \Delta_e, \Phi_e^+$,\dots introduced so far
really only depend on $J$, with the
exception of $\Phi_e^\circ$ which does involve $e$ itself: it is the
root system of the reductive Lie algebra $\c_e$. Let
$\Delta_e^\circ$ denote the base of $\Phi_e^\circ$ associated to the
positive system $\Phi_e^\circ \cap \Phi_e^+$. The {\em dominant
chamber} of the hyperplane arrangement $\mathscr A_e^\circ$
in $E_e$ associated to the root system $\Phi_e^\circ$
is $\{\alpha \in E_e \mid (\alpha,\beta)
> 0 \text{ for all } \beta \in \Delta_e^\circ\}$. It is usually {\em
not} a single chamber in the hyperplane arrangement $\mathscr A_e$
defined by $\Phi_e$, though it certainly contains the {\em standard
chamber} $\{\alpha \in E_e \mid (\alpha,\beta) > 0 \text{ for all }
\beta \in \Delta_e\}$ of $\mathscr A_e$ as a subset.

This setup gives a natural way to understand the
restricted Weyl group $W^J \cong N_W(W_J) / W_J$ from the previous section.
Recall that this acts faithfully on the vector space $E_e = E^J$,
so extending scalars we can view $W^J$ as a subgroup of $GL(\t_e)$. Let
$$
W_e = N_{G_e}(\t_e) / Z_{G_e}(\t_e),
$$
also naturally a subgroup of $GL(\t_e)$.
Using the decomposition $G_e = C_e \ltimes R_e$
and
noting that $\t_e \subseteq \c_e$,
it is easy to see that $N_{G_e}(\t_e) = N_{C_e}(\t_e) \ltimes Z_{R_e}(\t_e)$
and $Z_{G_e}(\t_e) = Z_{C_e}(\t_e) \ltimes Z_{R_e}(\t_e)$.
Hence, we can also write $W_e = N_{C_e}(\t_e) / Z_{C_e}(\t_e)$
as subgroups of $GL(\t_e)$.

\begin{lem}\label{lem1}
As subgroups of $GL(\t_e)$, we have that $W_e = W^J$.
\end{lem}

\begin{proof}
We identify $W$ with $N_G(T) / T = N_G(\t)/T$.

Take $x \in W^J$ represented by $\dot x \in N_G(T)$. Since $x\cdot
\Delta_J = \Delta_J$, $\dot x$ normalizes $L$.
Hence, $\dot x \cdot e$ is another
distinguished nilpotent element of $[\l,\l]$. We claim that there exists
$y \in L$ such that $\dot x \cdot e = y \cdot e$. To see this, it
suffices by the classification of distinguished nilpotent orbits in
$[\l,\l]$ to see that $\dot x \cdot e$ has the same labelled Dynkin
diagram as $e$. This is true because by inspection of the tables in
\cite{C}, the labelled Dynkin diagrams parameterizing distinguished
nilpotent orbits of Levi subalgebras of simple Lie algebras are
invariant under graph automorphisms. Hence, we have found an element
$y^{-1} \dot x \in G_e$ which normalizes $\t_e$ and acts on $\t_e$
in the same way as $x$. This shows that $W^J \subseteq W_e$.

Conversely, take $x \in W_e$ represented by $\dot x \in
N_{G_e}(\t_e)$. Recalling that $L = Z_G(\t_e)$, $\dot x$ certainly
normalizes $L$ too. Now $y \cdot T$ is a maximal torus of $L$, so
there exists
%It is well known that $N_G(L) = L(N_G(T) \cap N_G(L))$,
%so we can find some
$y \in L$ such that $y^{-1}\dot x  \in N_G(T)$.
Since $\dot x$ normalizes $\l$, it normalizes the center
$\t_e$ of $\l$, while $y$ centralizes $\t_e$.
Hence, $y^{-1}\dot x $ normalizes $\t_e$ and it acts on $\t_e$
in the same way as $x$. So $W_e \subseteq W^J$.
\end{proof}

Note that $W_e$ leaves $\Phi_e \subset E_e$ invariant, hence it acts
on bases for $\Phi_e$, or equivalently, on the chambers of the
hyperplane arrangement $\mathscr A_e$, as described by
Lemma~\ref{enum}. Let $W_e^\circ$ denote the Weyl group of the
reductive part $\c_e$ of $\g_e$, 
so $W_e^\circ$ is the subgroup of $GL(E_e)$ generated by the reflections in
the hyperplanes orthogonal to the simple roots $\Delta_e^\circ$ of
the root system $\Phi_e^\circ$ of $\c_e$. Let $Z_e$ denote the
stabilizer  in $W_e$ of the dominant chamber $\{\alpha \in E_e \mid
(\alpha,\beta) > 0 \text{ for all } \beta \in \Delta_e^\circ\}$.

\begin{lem}\label{lem11}
We have that $W_e = Z_e \ltimes W_e^\circ$
and $Z_e \cong C_e / C_e^\circ Z_{C_e}(\t_e)$,
a quotient of the component group $C_e / C_e^\circ \cong G_e / G_e^\circ$.
\end{lem}

\begin{proof}
Note that $W_e^\circ = N_{C_e^\circ}(\t_e) / Z_{C_e^\circ}(\t_e)
\cong N_{C_e^\circ}(\t_e) Z_{C_e}(\t_e) / Z_{C_e}(\t_e)$. Hence,
recalling that $W_e = N_{C_e}(\t_e) / Z_{C_e}(\t_e)$, we see that
the reflection group $W_e^\circ$ is a normal subgroup of $W_e$. Now
one can see that $W_e = Z_e\ltimes W_e^\circ$; see \cite[Lemma
2]{How}. Moreover, we have shown that $Z_e \cong W_e / W_e^\circ
\cong N_{C_e}(\t_e) / N_{C_e^\circ}(\t_e) Z_{C_e}(\t_e)$. Now
consider the natural map $N_{C_e}(\t_e) \rightarrow C_e / C_e^\circ
Z_{C_e}(\t_e)$. It is surjective because for every $x \in C_e$ there
exists $y \in C_e^\circ$ with $x \cdot \t_e = y \cdot \t_e$. Its
kernel is $N_{C_e^\circ}(\t_e) Z_{C_e}(\t_e)$. Hence it induces an
isomorphism between $Z_e$ and $C_e / C_e^\circ Z_{C_e}(\t_e)$.
\end{proof}

It follows from this and Lemma~\ref{enum} that $Z_e$ has $|\mathscr K_J|$
orbits on the set of chambers of the arrangement $\mathscr A_e$ that
are contained in the dominant chamber $\{\alpha \in E_e \mid
(\alpha,\beta) > 0 \text{ for all } \beta \in \Delta_e^\circ\}$, and
each orbit is regular. One can easily read off the structure of the
group $Z_e$ from the tables in \cite{How} and \cite{C}. For $\g$
simple, the group $Z_e$ is trivial, except in the following cases:
\begin{itemize}
\item[(i)] $\g = \sp_{2n}(\C)$ and $\lambda$ has $k > 0$
distinct even parts of even multiplicity, in which case
$Z_e \cong S_2 \times\cdots \times S_2\:(k\text{ times})$;
\item[(ii)] $\g = \so_N(\C)$, at least one part of $\lambda$ has odd multiplicity,
and $\lambda$ has $k > 0$
distinct odd parts of even multiplicity, in which case
$Z_e \cong S_2\times\cdots\times S_2 \:(k\text{ times})$;
\item[(iii)] $\g = \so_N(\C)$, all parts of $\lambda$ are of even multiplicity,
and $\lambda$ has $k > 1$
distinct odd parts,
in which case $Z_e \cong S_2\times\cdots\times S_2 \:((k-1)\text{ times})$;
\item[(iv)] $\g = F_4$ and $e$ has Bala-Carter label
$\widetilde A_1$, $A_2$ or $B_2$, in which case $Z_e \cong S_2$;
\item[(v)] $\g = E_6$ and
$e$ has Bala-Carter label $A_2$, in which case
$Z_e \cong S_2$.
\item[(vi)] $\g = E_6$ and
$e$ has Bala-Carter label $D_4(a_1)$, in which case
$Z_e \cong S_3$.
\item[(vii)] $\g = E_7$ and
$e$ has Bala-Carter label $A_2$, $A_2+A_1$, $D_4(a_1)+A_1$, $A_3+A_2$, $A_4$, $A_4+A_1$, $D_5(a_1)$ or $E_6(a_1)$, in which case
$Z_e \cong S_2$.
\item[(viii)] $\g = E_7$ and
$e$ has Bala-Carter label $D_4(a_1)$, in which case
$Z_e \cong S_3$.
\item[(ix)] $\g = E_8$ and $e$ has Bala-Carter label
$A_2$, $A_2+A_1$, $2A_2$,
$A_3+A_2$, $A_4$, $D_4(a_1)+A_2$, $A_4+A_1$, $D_5(a_1)$,
$A_4+2A_1$, $D_4+A_2$, $D_6(a_2)$,
$D_6(a_1)$, $E_6(a_1)$, $D_5+A_2$, $D_7(a_2)$, $E_6(a_1)+A_1$
or $D_7(a_1)$, in which case
$Z_e \cong S_2$.
\item[(x)] $\g = E_8$ and $e$ has Bala-Carter label
$D_4(a_1)$ or $D_4(a_1)+A_1$, in which case
$Z_e \cong S_3$.
\end{itemize}
In (i)--(iii), when $\g$ is classical, the partition $\lambda =
(1^{m_1} 2^{m_2} \cdots)$ denotes the Jordan type of $e$
in its natural representation.

Finally in this section, we wish to say a little more about the
dimensions of the {\em root spaces} $\g_e(\alpha,i)$ of $\g_e$. Note
that $\dim \g_e(\alpha,i)$ is the same as the multiplicity
$m(\alpha,i)$ from the proof of Lemma~\ref{wts}. By the definition
of $\c_e$, $m(\alpha,0)$ is one or zero according to whether $\alpha
\in \Phi_e^\circ$ or not. The root multiplicities $m(\alpha,i)$ for
$i > 0$ can often be greater than one, and in fact can be
arbitrarily large for symplectic and orthogonal Lie algebras. For
$\g = \sl_n$, the multiplicities $m(\alpha,i)$ are always one, and
explicit calculations as described in the next paragraph show that
$m(\alpha,i)$ is always at most three for $G$ simple of
exceptional type.
%it can only actually be equal to three if $G = E_8$ and $e$ has
%Bala-Carter label $E_6(a_3)+A_1$ or $E_7(a_5)$.
In general, the root space $\g_e(\alpha,i)$ need not be a subalgebra
of $\g_e$.

Let us explain exactly how to compute the root multiplicities
$m(\alpha,i)$ from the root system of $G$. Assume for this that $I =
\{1,\dots,r\} \setminus J$, so that as in the previous section
$\{\alpha_i^J \mid i \in I\}$ is the standard base for the
restricted root system $\Phi^J$. Of course, the restriction
$\alpha^J$ of a root $\alpha = \sum_{i=1}^r a_i \alpha_i$ is simply
$\sum_{i \in I} a_i \alpha_i^J$, so it is easy to write down the set
$\Phi^J$ explicitly given the root system of $\g$. Since $L$
centralizes $\t_e$, it does no harm to conjugate by an element of
$L$ to assume that the distinguished $\sl_2$-triple $(e,h,f)$ in
$[\l,\l]$ is in standard form, so that the values $\alpha_j(h)$ for
$j \in J$ are all either $0$ or $2$ as can be read off from the
labelled diagram for the distinguished nilpotent $e \in [\l,\l]$
from \cite{C}. Solving some linear equations, one can then uniquely
determine the other values $\alpha_i(h)$ for $i \in I$, using that
$h$ is orthogonal to $\t_e$. Hence we can compute all the integers
$\beta(h)$ for $\beta \in \Phi$. Now take $\alpha \in \Phi_e$. The
formal character of the $\s$-module arising from the $\alpha$-weight
space of $\g$ with respect to $\t_e$ is then
$$
\sum_{\beta \in \Phi\text{ s.t.\,}\beta^J = \alpha} x^{\beta(h)}.
$$
By $\sl_2$-theory, this can be written uniquely as $\sum_{i \geq 0}
m(\alpha,i) (x^i + x^{i-2}+\cdots+ x^{-i})$ for integers
$m(\alpha,i) \geq 0$. These are the desired multiplicities. For
exceptional groups, this procedure is particularly effective, and we
have implemented it in {\sc GAP} \cite{GAP4} to quickly compute all
root multiplicities in all cases, though there does not seem to be a
compact way to present this information here.
For classical groups, there is
a different approach,
described in sections \ref{sl_n}--\ref{so_N}.

\iffalse
\begin{example}\label{eg2}\rm
Take $G = E_7$ and $e$ with Bala-Carter
label $A_4+A_1$. Label the simple roots/the vertices of the Dynkin diagram as follows:
$$
\begin{array}{lllllll}
3&2&4&5&6&7\\
&&1
\end{array}
$$
Take $I = \{1,2\}$ and $J = \{3,4,5,6,7\}$.
Then,
$\Phi_e^+ = \{\alpha_1^J, \alpha_2^J, \alpha_1^J+\alpha_2^J,
\alpha_1^J+2\alpha_2^J, 2\alpha_1^J+2\alpha_2^J,
2\alpha_1^J+3\alpha_2^J\}$.
The values of $\alpha_i(h)$ for $i=1,\dots,7$ are given by the
labelled Dynkin diagram
$$
\begin{array}{rrrrrrr}
2&\!\!\!\,\text{-5}&\!2&2&2&2&2\\
&&\!\!\!\!\!\!\text{-4}
\end{array}
$$
From this, one computes the root multiplicities $m(\alpha,i)$ by the method explained above. To record these,
we list for
every $\alpha \in \Phi_e^+$ the sequence
made up of $m(\alpha,i)$ $i$'s for all $i \geq 0$:
$\alpha_1: 4$; $\alpha_2: 3,5$;
$\alpha_1+\alpha_2: 1,3,5,7$;
$\alpha_1+2\alpha_2: 2,6$;
$2\alpha_1+2\alpha_2: 4$;
$2\alpha_1+3\alpha_2: 1$.
$$
\begin{array}{|l|l|}
\hline
\alpha_1^J&4\\
\alpha_2^J&3,5\\
\alpha_1^J+\alpha_2^J&1,3,5,7\\
\alpha_1^J+2\alpha_2^J&2,6\\
2\alpha_1^J+2\alpha_2^J&4\\
2\alpha_1^J+3\alpha_2^J&1\\
\hline
\end{array}
$$
The restricted Cartan matrix with $ij$-entry
$\frac{2(\alpha_i^J, \alpha_j^J)}{(\alpha_j^J, \alpha_j^J)}$
for $i,j\in I$
is the matrix $\left(\begin{array}{rr}2 & -\frac{16}{7}\\ -\frac{4}{3}&2\;\end{array}\right)$.
\end{example}
\fi

\begin{example}\label{eg2}\rm
Take $G = E_7$ and $e$ with Bala-Carter
label $A_3+A_2$. Continuing with the notation from Example~\ref{eg1},
the values of $\alpha_i(h)$ for $i=1,\dots,7$ are given by the
labelled Dynkin diagram
$$
\begin{array}{rrrrrr}
2&2&\!\!\!\!\text{-5}&2&2&2\\
&&0
\end{array}
$$
From this, one computes the root multiplicities $m(\alpha,i)$ by the
method just explained. To record these, we list for every $\alpha
\in \Phi_e^+$ the sequence made up of $m(\alpha,i)$ $i$'s for all $i
\geq 0$: $\alpha^J_1:0$; $\alpha^J_2: 1,3,5$;
$\alpha^J_1+\alpha^J_2: 1,3,5$; $\alpha^J_1+2\alpha^J_2: 2,2,4,6$;
$\alpha^J_1+3\alpha^J_2: 3$; $2\alpha^J_1+3\alpha^J_2: 3$;
$2\alpha^J_1+4\alpha^J_2: 2$. The reductive part $C_e$ of the
centralizer is of type $A_1 + T_1$, and $\Delta_e^\circ =
\{\alpha_1^J\}$. Hence the fundamental chamber is the right half
plane in the diagram from Example~\ref{eg1}, and $W_e^\circ \cong
S_2$ is generated by the reflection in the vertical axis. The group
$Z_e \cong S_2$ is generated by the reflection in the horizontal
axis.
\end{example}

\section{Good gradings}

Continue with notation as in the previous section. In particular, we
have fixed bases $\Delta =\{\alpha_1,\dots,\alpha_r\}$ for $\Phi$
and $\Delta_e = \{\alpha_i^J \mid i \in I\}$ for $\Phi_e$. We often
now represent an element $c \in \t$ as a tuple $(c_1,\dots,c_r)$ of
complex numbers, where $c_i = \alpha_i(c)$. Of course, we think of
this as a labelling of the vertices of the Dynkin diagram. By a {\em
grading} of $\g$, we always mean an $\R$-grading
$$
\Gamma: \g = \bigoplus_{j \in \R} \g_j
$$
such that $[\g_i, \g_j] \subseteq \g_{i+j}$. We say that the grading
$\Gamma$ is {\em compatible with $\t$} if $\t \subseteq \g_0$. Since
every derivation of $\g$ is inner, there exists a unique semisimple
element $c \in \g$ defining $\Gamma$, i.e.\ so that $\g_j$ is the
$j$-eigenspace of $\ad c$. The grading is compatible with $\t$ if
and only if this element belongs to $\t$. In this way, gradings of
$\g$ that are compatible with $\t$ are parameterized by labelled
Dynkin diagrams $(c_1,\dots,c_r)$ with all labels $c_i \in \R$.
Every semisimple element of $\g$ is $G$-conjugate to an element of
$\t$, so every grading is $G$-conjugate to a grading that is
compatible with $\t$. Finally, two elements of $\t$ are
$G$-conjugate if and only if they are $W$-conjugate, hence every
grading of $\g$ is $G$-conjugate to a unique grading that is
compatible with $\t$ and whose labelled Dynkin diagram
$(c_1,\dots,c_r)$ has all labels $c_i \in \R_{\geq 0}$. We call this
labelled diagram the {\em characteristic} of the grading.

Now assume that $\Gamma:\g=\bigoplus_{j \in \R} \g_j$ is a good grading for
$e$ as defined in the introduction.
The proof of \cite[Theorem 1.3]{elakac} shows that
$\ad e: \g_j \rightarrow \g_{j+2}$ is surjective if and only if
$\ad e:\g_{-j-2} \rightarrow \g_{-j}$ is injective.
Hence the conditions that
$\ad e: \g_j \rightarrow \g_{j+2}$ is injective for all $j \leq -1$
and that
$\ad e: \g_j \rightarrow \g_{j+2}$ is surjective for all $j \geq -1$
in the definition of good grading are in fact equivalent.
So, $\Gamma$ is a good grading for $e \in \g_2$ if and only if
$\g_e \subseteq \bigoplus_{j > -1} \g_j$.

\begin{lem}\label{dimy}
Let $\Gamma$ be a grading of $\g$ with $e \in \g_2$. Then, $\dim
\g_e \geq \sum_{-1 \leq j < 1} \dim \g_j$ with equality if and only
if $\Gamma$ is a good grading for $e$.
\end{lem}

\begin{proof}
Copy the proof of \cite[Corollary 1.3]{elakac}.
\end{proof}

By a {\em good characteristic}, we mean the characteristic of a good
grading for $e$. By the proof of \cite[Theorem 1.2]{elakac}, a good
characteristic $(c_1,\dots,c_r)$ always has the property that $0
\leq c_i \leq 2$ for all $i=1,\dots,r$. We should observe that the
original good grading $\Gamma$ for $e$ can be recovered from its
characteristic uniquely up to conjugacy by $G_e$, i.e.\ good
characteristics parameterize $G_e$-conjugacy classes of good
gradings for $e$. To see this, suppose that $\Gamma$ and $\Gamma'$
are two good gradings for $e$ with the same characteristic. There
certainly exists $y \in G$ such that $y \cdot \Gamma' = \Gamma$. So
$\Gamma$ is good both for $e$ and for $y \cdot e$. Letting $G_0$ be
the set of all elements of $G$ that preserve the grading $\Gamma$,
Lemma~\ref{dense} below implies that $y \cdot e = z \cdot e$ for
some $z \in G_0$. But then $z^{-1} y \cdot \Gamma' = \Gamma$ too,
and $z^{-1} y \in G_e$, as required.

\begin{lem}\label{dense}
If $\Gamma$ is a good grading,
the set of all elements $e \in \g_2$
such that $\Gamma$ is a good grading for $e$ is a
dense open orbit for the action of $G_0$ on $\g_2$.
\end{lem}

\begin{proof}
Suppose that $\Gamma$ is a good grading for $e$ and for $e'$.
We have that $[e, \g_0] = \g_2$.
Hence $\dim G_0 \cdot e = \dim \g_2$, and $G_0 \cdot e$ is dense open
in $\g_2$. So is $G_0 \cdot e'$, so $G_0 \cdot e$ and $G_0 \cdot e'$ have non-empty intersection. Hence,
$G_0 \cdot e = G_0 \cdot e'$.
\end{proof}

For the next lemma, we note that $E_e$ is the
$\R$-form for $\t_e$ consisting of all $p \in \t_e$
such that the eigenvalues of $\ad p$ on $\g$ are real.
Also recall that the element $h$ from our
fixed $\sl_2$-triple $(e,h,f)$
belongs to $\t$, since $\t$ was chosen originally
to lie in $\c$.

\begin{lem}\label{choice}
Every $G_e$-conjugacy class of good gradings for $e$
has a representative $\Gamma$ that is compatible with $\t$.
Moreover, for any such $\Gamma$,
we have that $h \in \g_0, f \in \g_{-2}$, and the element $c \in \g$
defining the grading $\Gamma$ is of the form $c = h + p$ for some
point $p \in E_e$.
\end{lem}

\begin{proof}
Let $\Gamma$ be any good grading for $e$.
As in \cite[Lemma 1.1]{elakac},
there
exists an $\mathfrak{sl}_2$-triple $(e,h',f')$ with $h' \in \g_0$ and
$f' \in \g_{-2}$.
This is conjugate to our fixed $\mathfrak{sl}_2$-triple
$\s = (e,h,f)$ by an element $x$ of $G^\circ_e$.
On replacing $\Gamma$ by $x \cdot \Gamma$ if necessary, we
may therefore assume that $h \in \g_0$ and $f \in \g_{-2}$ already.
Let $c$ be
the semisimple element of $\g$ defining the grading $\Gamma$.
Since $[h,e]=[c,e]=2e$, the element $p=c-h$
centralizes $e$.
Since $h \in \g_0$, we have that $[h,c] = 0$,
so $p$ is a semisimple element of $\c_e$.
Recalling that $\t_e$ is a Cartan subalgebra of $\c_e$, we can
therefore conjugate once more by an element of $C^\circ_e$ to
reduce to the situation that $p \in \t_e$.
But then $c = h+p$ belongs to $\t$, and the grading
$\Gamma$ is compatible with $\t$ as required.

Now suppose that $\Gamma$ is any good grading for $e$ that is
compatible with $\t$. Let $c$ be the element of $\t$ defining the
grading, so that $[c,e]=2e$ and $[c,h] = 0$, i.e.\ $h \in \g_0$.
This implies that $p=c-h$ centralizes $e$ and $h$, hence also $f$,
so $[c,f]=[h,f]=-2f$ and $f \in \g_{-2}$. Finally, observe that $p$
belongs to $\t_e$, hence to $E_e$ since $\Gamma$ is an $\R$-grading.
\end{proof}

Given any $p \in E_e$, let $\Gamma(p)$ denote the grading of $\g$
defined by the semisimple element $h+p$. For example, $\Gamma(0)$ is
the Dynkin grading. In general, the grading $\Gamma(p)$ is certainly
compatible with $\t$ and the element $e$ is in degree $2$. However,
it need not be a good grading for $e$.

\begin{thm}\label{ggp}
For $p \in E_e$, the grading $\Gamma(p)$ is a good grading for $e$
if and only if $|\alpha(p)| < d(\alpha)$ for all $\alpha \in
\Phi^+_e$, where $d(\alpha) = 1+\min\{i \geq 0 \mid m(\alpha,i) \neq
0\}$, that is, the minimal dimension of an irreducible
$\s$-submodule of the $\alpha$-weight space of $\g$ with respect to
$\t_e$.
\end{thm}

\begin{proof}
Recall that $\Gamma(p):\g = \bigoplus_{j \in \Z} \g_j$ is a good
grading for $e$ if and only if $\g_e \subseteq \bigoplus_{j > -1}
\g_j$. Let $\g = \bigoplus_{\alpha \in \Phi_e\cup\{0\}} \bigoplus_{i
\geq 0} m(\alpha,i) L(\alpha,i)$ be the decomposition of $\g$ as an
$\t_e \oplus \s$-module, as in the proof of Lemma~\ref{wts}. Note
that $h+p$ acts on the highest weight vector of $L(\alpha,i)$ as the
scalar $\alpha(p)+i$. Putting these things together, we get that
$\Gamma(p)$ is a good grading for $e$ if and only if $\alpha(p)+i >
-1$ whenever $m(\alpha,i) \neq 0$. Since $\Phi_e = \Phi_e^+ \sqcup
(-\Phi_e^+)$, the theorem follows easily.
\end{proof}

Let $\mathscr P_e$ denote the set of all $p \in E_e$ such that
$\Gamma(p)$ is a good grading for $e$.
Since $\Phi^+_e$ spans $E_e^*$,
Theorem~\ref{ggp}
shows in particular that $\mathscr P_e$ is an open
convex polytope in the real vector space
$E_e$. It can be computed explicitly from information about the
root multiplicities $m(\alpha,i)$; see the discussion before
Example~\ref{eg2}.
We call it the
{\em good grading polytope} corresponding to $e$.
By Lemma~\ref{choice}, the map $p \mapsto \Gamma(p)$ gives a
bijection between $\mathscr P_e$
and the set of
all good gradings for $e$ that are compatible with
$\t$.
The next theorem describes exactly when the good gradings $\Gamma(p)$
and $\Gamma(p')$ for $e$ are conjugate,
for points $p, p' \in \mathscr P_e$. Recall this is if and only if
they have the same characteristic.

\begin{thm}\label{conj}
For $p, p' \in \mathscr P_e$,
the good gradings $\Gamma(p)$ and $\Gamma(p')$ are
$G$-conjugate if and only if $p$ and $p'$ are $W_e$-conjugate.
\end{thm}

\begin{proof}
Suppose that the good gradings $\Gamma(p): \g = \bigoplus_{i \in \R}
\g_i$ and $\Gamma(p'):\g = \bigoplus_{j \in \R} \g_j'$ are
$G$-conjugate. Since they are both good gradings for $e$, they are
already conjugate under the centralizer $G_e$, as we explained
earlier. So we can find $x \in G_e$ such that $x\cdot (h+p) = h+p'$.
Since $h+p'$ lies in both $\t$ and in $x \cdot \t$, it centralizes
$\t_e$ and $x \cdot \t_e$, so both $\t_e$ and $x \cdot \t_e$ lie in
$\g_0'$. They also clearly both lie in $\g_e$, hence $\t_e$ and $x
\cdot \t_e$ are Cartan subalgebras of $\g_e\cap\g_0'$. Letting
$G_0'$ be the subgroup of $G$ consisting of all elements that
preserve the grading $\Gamma'$, we deduce that there exists an
element $y \in G_e \cap G_0'$ such that $y \cdot \t_e = x \cdot
\t_e$. Hence, $y^{-1} x \cdot (h+p) = h+p'$ and $y^{-1}x \in
N_{G_e}(\t_e)$. Since $y^{-1} x$ normalizes $\t_e$, it normalizes
$\l$, hence $[\l,\l]$. So, $y^{-1} x \cdot h = h'$ for some $h' \in
[\l,\l]$, and $h+p'=y^{-1}x \cdot (h+p) = h'+y^{-1}x \cdot p$. This
shows that $h=h'$ and $y^{-1}x\cdot p = p'$ already. Hence $p$ and
$p'$ are $W_e$-conjugate.

Conversely, suppose that $p$ and $p'$ are $W_e$-conjugate. Then,
recalling that $W_e = N_{C_e}(\t_e) / Z_{C_e}(\t_e)$, we can find $x
\in N_{C_e}(\t_e)$ with $x \cdot p = p'$. Since $x$ lies in $C_e$ it
centralizes $h$, so $x \cdot (h+p) = h+p'$. Hence $\Gamma(p)$ and
$\Gamma(p')$ are conjugate.
\end{proof}

We finally introduce the {\em affine hyperplanes}
$$
H_{\alpha,k} = \{p \in E_e \mid \alpha(p) = k\}
$$
for each $\alpha \in \Phi_e^+$ and $k \in \Z$. The significance of
these will be discussed in detail in the next section. We just want
to point out here that the {\em integral} good gradings for $e$ that
are compatible with $\t$ are parameterized by the points $p \in
\mathscr P_e$ such that $\alpha(p) \in \Z$ for all $\alpha \in
\Phi_e^+$. In other words, $\Gamma(p)$ is an integral grading if and
only if $p$ lies  on the same number of the affine hyperplanes
$H_{\alpha,k}$ for $\alpha \in \Phi_e^+,k\in\Z$ as the origin (which
corresponds to the Dynkin grading). Actually, it is often the case
that the Dynkin grading is the only integral good grading, as is
well explained by \cite[Corollary 1.1]{elakac}.

\begin{example}\label{eg3}\rm
Continue with $G = E_7$ and $e$ having Bala-Carter label $A_3+A_2$,
notation as in Examples~\ref{eg1} and \ref{eg2}.
In Example~\ref{eg2},
we computed {all} the root multiplicities $m(\alpha,i)$.
Hence, according to Theorem~\ref{ggp}, the good grading polytope
$\mathscr P_e$ is the subspace of $E_e$ defined by the inequalities
\begin{align*}
|\alpha_1^J(p)| < 1,\quad |\alpha_2^J(p)| < 2, \quad
|(\alpha_1^J+\alpha_2^J)(p)| < 2,\quad &|(\alpha_1^J+2\alpha_2^J)(p)| < {3},\\
|(\alpha_1^J  + 3 \alpha_2^J)(p)| < 4,\quad |(2
\alpha_1^J+3\alpha_2^J)(p)| < 4,\quad
&|(2\alpha_1^J+4\alpha_2^J)(p)| < 3.
\end{align*}
These are equivalent just to the inequalities $|\alpha_1^J(p)| < 1,
|(\alpha_1^J+2\alpha_2^J)(p)| < \frac{3}{2}$, so the good grading
polytope (drawn on the same axes as in Example~\ref{eg1} but with a
different scale) can be represented as the interior of the
rectangle:
$$
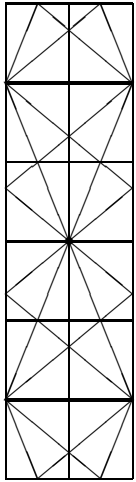
\begin{picture}(150,182.25)
\put(75,89.75){\line(1,0){24}}
\put(75,89.75){\line(-1,0){24}}
\put(75,90){\line(1,0){24}}
\put(75,90){\line(-1,0){24}}
\put(75,90.25){\line(1,0){24}}
\put(75,90.25){\line(-1,0){24}}
\put(75,120){\line(1,0){24}}
\put(75,120){\line(-1,0){24}}
\put(75,150){\line(1,0){24}}
\put(75,150){\line(-1,0){24}}
\put(75,149.75){\line(1,0){24}}
\put(75,149.75){\line(-1,0){24}}
\put(75,150.25){\line(1,0){24}}
\put(75,150.25){\line(-1,0){24}}
\put(75,180){\line(1,0){24}}
\put(75,180){\line(-1,0){24}}
\put(75,60){\line(1,0){24}}
\put(75,60){\line(-1,0){24}}
\put(75,30){\line(1,0){24}}
\put(75,30){\line(-1,0){24}}
\put(75,29.75){\line(1,0){24}}
\put(75,29.75){\line(-1,0){24}}
\put(75,30.25){\line(1,0){24}}
\put(75,30.25){\line(-1,0){24}}
\put(75,0){\line(1,0){24}}
\put(75,0){\line(-1,0){24}}
\put(75,90){\line(0,1){90}}
\put(75,90){\line(0,-1){90}}
\put(99,90){\line(0,1){90}}
\put(99,90){\line(0,-1){90}}
\put(51,90){\line(0,1){90}}
\put(51,90){\line(0,-1){90}}
\put(75,90){\line(-2,5){24}}
\put(75,90){\line(2,-5){24}}
\put(75,90){\line(2,5){24}}
\put(75,90){\line(-2,-5){24}}
\put(51,149.625){\line(2,5){12.15}}
\put(99,149.625){\line(-2,5){12.15}}
\put(51,30.375){\line(2,-5){12.15}}
\put(99,30.375){\line(-2,-5){12.15}}
\put(75,90){\line(-6,5){24}}
\put(75,90){\line(6,5){24}}
\put(75,90){\line(6,-5){24}}
\put(75,90){\line(-6,-5){24}}
\put(51.1,109.96875){\line(6,5){48.3}}
\put(51.1,70.18125){\line(6,-5){48.3}}
\put(98.9,109.96875){\line(-6,5){48.3}}
\put(98.9,70.18125){\line(-6,-5){48.3}}
\put(51.1,150.13125){\line(6,5){35.55}}
\put(98.9,150.13125){\line(-6,5){35.55}}
\put(51.1,29.90625){\line(6,-5){35.55}}
\put(98.9,29.90625){\line(-6,-5){35.55}}
\put(75,90){\circle*{3}}
\end{picture}
$$
We have also drawn on the diagram
the affine hyperplanes
that intersect with $\mathscr P_e$. From this, we see
that no other point of
$\mathscr P_e$ lies on as many affine hyperplanes as the origin.
So the only integral good grading for $e$ is the Dynkin grading,
which is consistent with the classification of
integral good gradings from \cite{elakac} in this case.
\end{example}

\section{Alcoves and adjacencies}

Recall from the introduction that a pair $\Gamma: \g = \bigoplus_{i
\in \Z} \g_i$ and $\Gamma': \g = \bigoplus_{j \in \Z} \g_j'$ of good
gradings for $e$ are {\em adjacent} if $\g = \bigoplus_{i^- \leq j
\leq i^+} \g_i \cap \g_j'$. Note this means in particular that the
gradings $\Gamma$ and $\Gamma'$ are compatible with each other,
i.e.\ the semisimple elements of $\g$ that define the gradings
$\Gamma$ and $\Gamma'$ commute. Moreover, $\bigoplus_{i, j \geq 0}
\g_i \cap \g_j'$ is a parabolic subalgebra of $\g$ with Levi factor
$\g_0 \cap \g_0'$. So we can find an element $x \in G$ such that $x
\cdot \t \subseteq \g_0 \cap \g_0'$ and $x\cdot \b \subseteq
\bigoplus_{i, j \geq 0} \g_i \cap \g_j'$. The characteristics
$(c_1,\dots,c_r)$ and $(c_1',\dots,c_r')$ of the gradings $\Gamma$
and $\Gamma'$ can then be read off simultaneously from the equations
$x \cdot \g_{\alpha_i} \subseteq \g_{c_i} \cap \g_{c_i'}$ for all
$i=1,\dots,r$.

We say that two good characteristics $(c_1,\dots,c_r)$ and
$(c_1',\dots,c_r')$ are {\em adjacent} if they are the characteristics
of a pair of adjacent good gradings $\Gamma$ and $\Gamma'$ for $e$.
If we are given just a pair $(c_1,\dots,c_r)$ and
$(c_1',\dots,c_r')$ of adjacent good characteristics, we can recover the pair
$\Gamma$ and $\Gamma'$ of adjacent good gradings for $e$
from which they were defined uniquely up to
conjugation by $G_e$. To see this, the previous paragraph implies that
$\Gamma$ and $\Gamma'$ can be obtained by simultaneously
conjugating the two gradings defined by
declaring that each $\g_{\pm \alpha_i}$ is in degree $\pm c_i$ or
in degree $\pm c_i'$, respectively, by some element $x \in G$.
So we just need to show that if $\Lambda$ and $\Lambda'$ are another
pair of good gradings for $e$ with $y \cdot \Lambda = \Gamma$
and $y \cdot \Lambda' = \Gamma'$ for some $y \in G$,
then in fact $\Lambda$ and $\Lambda'$ are already
conjugate to $\Gamma$ and $\Gamma'$
by an element of $G_e$.
For this, note that $\Gamma$ and $\Gamma'$ are good gradings both for
$e$ and for $y \cdot e$. So by Lemma~\ref{ano} below, there exists an
element $z \in G$ preserving both gradings $\Gamma$ and $\Gamma'$ with
$z \cdot e = y \cdot e$. But then
$z^{-1} y \cdot \Lambda = \Gamma$ and
$z^{-1} y \cdot \Lambda' = \Gamma'$, and we have that
$z^{-1} y \in G_e$, as required.

\begin{lem}\label{ano}
Let $\Gamma$ and $\Gamma'$ be adjacent good gradings for $e$.
Let $G_0$ and $G_0'$ be the subgroups of $G$ consisting of
all elements that preserve the gradings
$\Gamma$ and $\Gamma'$, respectively.
The set of all elements $e'$ of $\g_2 \cap \g_2'$ such that both $\Gamma$
and $\Gamma'$ are good gradings for $e'$ is a dense open orbit
for the action of the group $G_0\cap G_0'$.
\end{lem}

\begin{proof}
Let $\g_{< j}$ denote $\bigoplus_{i < j} \g_i$.
Define $\g_{> j}, \g'_{< j}$ and $\g'_{> j}$ similarly.
Since the map $\ad e: \g_0 \rightarrow \g_2$ is surjective and
preserves the direct sum decompositions $\g_0 = (\g_0 \cap \g_{<0}')
\oplus (\g_0 \cap \g_0') \oplus (\g_0\cap\g_{>0}')$ and $\g_2 =
(\g_2 \cap \g_{<2}') \oplus (\g_2 \cap \g_2') \oplus
(\g_2\cap\g_{>2}')$, we have that $[\g_0 \cap \g_0', e] = \g_2 \cap
\g_2'$. Now argue as in the proof of Lemma~\ref{dense}.
\end{proof}

This shows that adjacent good characteristics parameterize
$G_e$-conjugacy classes of adjacent good gradings for $e$. In order
to classify all adjacent good gradings, hence all adjacent good
characteristics, in terms of the good grading polytope we use the
following lemma.

\begin{lem}\label{wma}
Let $\Gamma$ and $\Gamma'$
be adjacent good gradings for $e$.
Then there exists $x \in G_e^\circ$ such that both $x \cdot\Gamma$
and $x \cdot \Gamma'$ are compatible with $\t$.
\end{lem}

\begin{proof}
Note $e \in \g_2 \cap \g_2'$ by the definition of good grading.
Arguing as in the proof of \cite[Lemma 1.1]{elakac}, working with
the bigrading $\g = \bigoplus_{i,j \in \Z} \g_i \cap \g_j'$ instead
of the grading used there, one shows that there exists an
$\mathfrak{sl}_2$-triple $(e, h', f')$ with $h' \in \g_0 \cap \g_0'$
and $f' \in \g_{-2} \cap \g_{-2}'$. Conjugating by an element of
$G_e^\circ$ if necessary, we may assume that in fact $h'=h$ and
$f'=f$, i.e.\ $h \in \g_0 \cap \g_0'$ and $f \in \g_{-2}\cap
\g_{-2}'$ already. Now argue as in the proof of Lemma~\ref{choice}
to show that $\t \subseteq \g_0 \cap \g_0'$.
\end{proof}

Assume that $\Gamma$ and $\Gamma'$ are adjacent good gradings for
$e$. In view of Lemma~\ref{wma}, we can conjugate by an element of
$G_e^\circ$ if necessary to assume that $\Gamma$ and $\Gamma'$ are
both compatible with $\t$. Then, $\Gamma = \Gamma(p)$ and $\Gamma' =
\Gamma(p')$ for points $p,p' \in \mathscr P_e$. In this way, the
problem of determining all conjugacy classes of pairs of adjacent
good gradings for $e$ reduces to describing exactly when $\Gamma(p)$
and $\Gamma(p')$ are adjacent. To do this, recall the affine
hyperplanes $\{H_{\alpha,k} \mid \alpha \in \Phi_e^+, k \in \Z\}$
introduced at the end of the previous section. We refer to the
connected components of $$ E_e \setminus \displaystyle
\bigcup_{\substack{\alpha \in \Phi^+_e\\ k \in \Z\:\:\:}}
H_{\alpha,k}
$$
as {\em alcoves}.
Note that the closure of
$\mathscr P_e$ is the union of the closures of finitely many alcoves.

\begin{thm}\label{adj} Let $p,p' \in \mathscr P_e$.
Then, $\Gamma(p)$ and $\Gamma(p')$ are adjacent if and only if $p$ and $p'$
belong to the closure of the same alcove.
\end{thm}

\begin{proof}
By definition, $\g_\alpha$ lies in the degree $\alpha(h)+\alpha(p)$
piece of the grading $\Gamma(p)$,
for each $\alpha \in \Phi$. Since $\alpha(h)$ is always an integer,
it follows that $\Gamma(p)$ and $\Gamma(p')$ are adjacent if and only
if $\alpha(p)^- \leq \alpha(p') \leq \alpha(p)^+$ for all $\alpha \in \Phi$.
Equivalently, $p$ and $p'$ belong to the closure of the same alcove.
\end{proof}

Theorems~\ref{ggp} and \ref{adj} combine to prove Theorem~\ref{thm2}
from the introduction. In the remainder of the section, we want to
explain the remaining ingredients needed to deduce
Theorem~\ref{thm1} from it, as we outlined in the introduction. The
first step is accomplished by the following lemma. Recall from the
introduction that $\langle.\,,.\rangle$ is the skew-symmetric
bilinear form defined from $\langle x,y\rangle = ([x,y],e)$.

\begin{lem}\label{lem3}
Let $\Gamma:\g=\bigoplus_{i \in \R} \g_i$ and
$\Gamma':\g=\bigoplus_{j \in \R} \g_j$ be adjacent good gradings for
$e$. Then there exist Lagrangian subspaces $\k$ of $\g_{-1}$ and
$\k'$ of $\g_{-1}'$ (both with respect to the form
$\langle.\,,.\rangle$) such that
$$
\k \oplus \bigoplus_{i < -1} \g_i =
\k' \oplus\bigoplus_{j < -1} \g_j'.
$$
\end{lem}

\begin{proof}
In the notation from the proof of Lemma~\ref{ano}, we have that
\begin{align*}
\g_{-1} &= ((\g_{-1}\cap \g_{<-1}')\oplus (\g_{-1} \cap \g_{>-1}'))
\perp (\g_{-1} \cap \g_{-1}'),\\
\g_{-1}' &= ((\g_{<-1}\cap \g_{-1}') \oplus (\g_{>-1} \cap \g_{-1}'))
\perp
(\g_{-1} \cap \g_{-1}'),
\end{align*}
where the $\perp$'s are with respect to the form
$\langle.\,,.\rangle$. Hence the restriction of
$\langle.\,,.\rangle$ to $\g_{-1} \cap \g_{-1}'$ is non-degenerate.
Let $\k''$ be a Lagrangian subspace of $\g_{-1} \cap \g_{-1}'$. Then
set $\k = \k'' \oplus (\g_{-1} \cap \g_{<-1}')$ and $\k' = \k''
\oplus (\g_{<-1} \cap \g_{-1}')$. This does the job in view of the
definition of adjacency.
\end{proof}

It just remains to indicate
how to adapt the argument of Gan and Ginzburg \cite{gangin} to our slightly more general setting.
For an integer $d \geq 1$, we define a {\em $d$-good grading} for $e$
to be a $\Z$-grading $\g = \bigoplus_{j \in \Z} \g_j$ with $e \in \g_d$, such
that $\ad e:\g_j \rightarrow \g_{j+d}$ is injective for
$j \leq -\frac{d}{2}$ and surjective
for $j \geq -\frac{d}{2}$.
Given any good $\R$-grading $\Gamma$ for $e$ in our old sense,
it is easy to see that one can always replace $\Gamma$ by a good $\Q$-grading
without changing $e \in \g_2$ or
either of the spaces $\g_{-1}$ or $\bigoplus_{j < -1} \g_j$,
which is all that matters for the construction of finite $W$-algebras.
In turn, since $\g$ is finite dimensional, any good $\Q$-grading
for $e \in \g_2$ can be scaled by a sufficiently large integer $d$
so that it becomes a $2d$-good
grading for $e \in \g_{2d}$ in the new sense.
This reduces to the situation that $\Gamma:\g = \bigoplus_{j \in \Z} \g_j$
is a $2d$-good grading for $e$.

Next let $\k$ be any {\em isotropic} subspace of $\g_{-d}$, and let
$\k^\perp \subseteq \g_{-d}$ be its annihilator with respect to the
form $\langle.\,,.\rangle$. Let $\m = \k \oplus \bigoplus_{i < -d}
\g_i$ and $\n = \k^\perp \oplus \bigoplus_{j < -d} \g_j$. Note
$\chi:\m \rightarrow \C, x \mapsto (x,e)$ is a representation of
$\m$. Set $Q_\k = U(\g) \otimes_{U(\m)} \C_\chi = U(\g) / I_\k$,
where $I_\k$ is the left ideal of $U(\g)$ generated by all
$\{x-\chi(x) \mid x \in \m\}$. Since $I_\k$ is stable under the
adjoint action of $\n$, we get an induced action of $\n$ on $Q_\k$,
by $x \cdot (u \otimes 1_\chi) = [x,u] \otimes 1_\chi$ for $x \in
\n, u \in U(\g)$. Let $H_\k$ be the space $Q_\k^{\n}$ of $\n$-fixed
points with respect to this action. It has a well-defined algebra
structure defined by $(u \otimes 1_\chi) (v \otimes 1_\chi) = (uv)
\otimes 1_\chi$, for $u,v \in U(\g)$ such that $u \otimes 1_\chi, v
\otimes 1_\chi \in Q_\k^{\n}$. Now we can formulate the slight
generalization of Gan and Ginzburg's theorem that we need here.

\begin{thm}\label{ggt}
Let $\k \subseteq \k'$ be two isotropic subspaces of $\g_{-d}$,
and define the corresponding algebras $H_\k$
and $H_{\k'}$ as above.
The natural map $Q_\k \rightarrow Q_{\k'}$ induced by
the inclusion $\k \hookrightarrow \k'$
restricts to an algebra isomorphism $H_\k \rightarrow H_{\k'}$.
\end{thm}

If we take $\k = \{0\}$ and $\k', \k''$ to be any two
Lagrangian subspaces of $\g_{-d}$, this theorem gives
isomorphisms from $H_\k$ to both $H_{\k'}$ and $H_{\k''}$.
Composing one with the inverse of the other, we get
a canonical isomorphism between $H_{\k'}$ and $H_{\k''}$.
In turn, by the Frobenius reciprocity argument
explained in the introduction of \cite{brukle},
$H_{\k'}$ is naturally isomorphic to the finite $W$-algebra
$H_{\chi'}$ from the introduction defined from the Lagrangian subspace
$\k'$, while $H_{\k''}$ is naturally isomorphic to $H_{\chi''}$
defined from $\k''$.
In this way, we obtain the canonical isomorphism
between $H_{\chi'}$ and $H_{\chi''}$ needed to prove Theorem~\ref{thm1}.

The proof of Theorem~\ref{ggt} itself is almost exactly the same as
the proof of the second part of \cite[Theorem 4.1]{gangin}. We just
note here that one needs to replace the linear action $\rho$ of
$\C^\times$ on $\g$ from \cite{gangin} with one defined by
$\rho(t)(x) = t^{d-j} x$, for $t \in \C^\times$ and $x \in \g_j$.
There is a corresponding Kazhdan filtration on $U(\g)$ as in
\cite[\S 4]{gangin}. We also note the identity
$$
\m^\perp = [\n,e] \oplus \z_\g(f),
$$
where $\m^\perp$ is annihilator of $\m$ in $\g$
with respect to the Killing form and $(e,h,f)$ is
an $\sl_2$-triple with $f \in \g_{-d}$; this is proved as in
\cite[(2.2)]{gangin} using Lemma~\ref{dimy}.
Combining these things, one gets the analogue of
\cite[Lemma 2.1]{gangin}, which is the
key lemma needed
in
the spectral sequence argument used to prove \cite[Theorem 4.1]{gangin}.

\section{Good gradings for $\sl_n(\C)$} \label{sl_n}

In this section, we describe explicitly the restricted root systems
and the good grading polytopes
arising from $\g = \sl_n(\C)$.
Let $V$ denote the natural $n$-dimensional $\g$-module of column
vectors,
with standard basis $v_1,\dots,v_n$. Also let $\t$ be the
standard Cartan subalgebra consisting of all diagonal matrices in $\g$.
Letting $\delta_i$ be the element of $\t^*$ picking out the $i$th
diagonal entry of a matrix in $\t$, the elements
$\alpha_1 = \delta_1-\delta_2,\dots,\alpha_{n-1} = \delta_{n-1}-\delta_n$
give a base $\Delta$ for the root system $\Phi$.

Nilpotent orbits in $\g$ are parameterized by partitions $\lambda =
(\lambda_1\geq\lambda_2 \geq \cdots)$ of $n$. Fix such a partition
$\lambda$ throughout the section having a total of $m$ non-zero
parts. In order to write down an $\sl_2$-triple corresponding to the
partition $\lambda$ explicitly, we first recall the definition of
the {\em Dynkin pyramid} of shape $\lambda$, following
\cite{elakac}. This is a diagram consisting of $n$ boxes each of
size $2$ units by $2$ units drawn in the upper half of the
$xy$-plane. By the coordinates of a box, we mean the coordinates of
its midpoint. We will also talk about the {\em row number} of a box,
meaning its $y$-coordinate, and the {\em column number} of a box,
meaning its $x$-coordinate. Letting $r_i = 2i-1$ for short, the
Dynkin pyramid has $\lambda_i$ boxes in row $r_i$ centered in
columns $1-\lambda_i,3-\lambda_i,\dots,\lambda_i-1$, for each
$i=1,\dots,m$. For example, here is the Dynkin pyramid of shape
$\lambda = (3,3,2)$:
$$
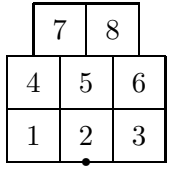
\begin{picture}(60,60)
\put(0,0){\line(1,0){60}}
\put(0,20){\line(1,0){60}}
\put(0,40){\line(1,0){60}}
\put(10,60){\line(1,0){40}}
\put(60,0){\line(0,1){40}}
\put(0,0){\line(0,1){40}}
\put(20,0){\line(0,1){40}}
\put(40,0){\line(0,1){40}}
\put(10,40){\line(0,1){20}}
\put(30,40){\line(0,1){20}}
\put(50,40){\line(0,1){20}}
\put(30,0){\circle*{3}}
\put(10,10){\makebox(0,0){{1}}}
\put(30,10){\makebox(0,0){{2}}}
\put(50,10){\makebox(0,0){{3}}}
\put(10,30){\makebox(0,0){{4}}}
\put(30,30){\makebox(0,0){{5}}}
\put(50,30){\makebox(0,0){{6}}}
\put(20,50){\makebox(0,0){{7}}}
\put(40,50){\makebox(0,0){{8}}}
\end{picture}
$$
We fix once and for all a numbering $1,2,\dots,n$ of the boxes of
the Dynkin pyramid, and let $\row(i)$ and $\col(i)$ denote the row
and column numbers of the $i$th box. Writing $e_{i,j}$ for the
$ij$-matrix unit, let $e = \sum_{i,j} e_{i,j}$ summing over all $1
\leq i,j \leq n$ such that $\row(i) = \row(j)$ and $\col(i) =
\col(j)+2$. This is a nilpotent matrix of Jordan type $\lambda$. For
example, taking $\lambda = (3,3,2)$ and numbering boxes as above, we
have that $e = e_{8,7}+e_{6,5}+e_{5,4}+e_{3,2}+e_{2,1}$. Also let $h
= \sum_{i=1}^n \col(i) e_{i,i}$. There is then a unique element $f
\in \g$ such that $(e,h,f)$ is an $\sl_2$-triple.

With these choices, it is the case that
$\t$ is contained in $\c$
and $\t_e$ is a Cartan subalgebra of $\c_e$,
as was required in section \ref{s3}.
Explicitly, $\t_e$ consists of all matrices in $\t$ such that
the $i$th and $j$th diagonal entries are equal
whenever $\row(i) = \row(j)$, and the real vector space $E_e$
consists of all such matrices with entries in $\R$.
Let $\eps_i \in E_e^*$ be the function picking out the
$j$th diagonal entry of a matrix in $E_e$, where
$j$ here is chosen so that the $j$th box of the Dynkin pyramid
is in row $r_i$.
Then, $E_e^*$ is the $(m-1)$-dimensional real vector space
spanned by $\eps_1,\dots,\eps_m$ subject to the relation
$\sum_{i=1}^m \lambda_i \eps_i = 0$.
It is natural to identify $E_e$ and $E_e^*$ via the real inner product arising
from the trace form on $\g$, with respect to which
$(\eps_i,\eps_j) = \frac{\delta_{i,j}}{\lambda_i}$.
We have that
$$
\Phi_e = \{\eps_i - \eps_j \mid 1 \leq i, j \leq m, i \neq j\}.
$$
Setting $p_i = \eps_i(p)$,
any point $p \in E_e$ can be represented
as a tuple $(p_1,\dots,p_m) \in \R^m$ with $\sum_{i=1}^m \lambda_i p_i = 0$.
The values of $d(\alpha)$ from Theorem \ref{ggp} can be determined
using an explicit description of $\g_e$ as in \cite{elakac}: we have that
$$
d(\eps_i - \eps_j) = 1 + |\lambda_i - \lambda_j|
$$
for all $1 \leq i, j \leq m$.  Therefore, the good grading polytope
$\mathscr P_e$ is the open subset of $E_e$ consisting of all points
$p=(p_1,\dots,p_m)$ such that
$$
|p_i - p_j| < 1 + \lambda_i - \lambda_j
$$
for all $1 \leq i < j \leq m$.
Also, the restricted Weyl group $W_e$ is the group $S_{m_1} \times S_{m_2}
\times\cdots$, where $m_i$ denotes the number of parts of $\lambda$
that equal $i$, acting on $E_e$ by permuting all $\eps_i$'s of equal
length.

As explained in \cite{elakac}, there is a convenient way to
visualize the good grading $\Gamma(p)$ corresponding to
$p=(p_1,\dots,p_m) \in \mathscr P_e$. First, associate a pyramid
$\pi(p)$ to $p$ by sliding all numbered boxes in row $r_i$ of the
Dynkin pyramid to the right by $p_i$ units, for each $i=1,\dots,m$;
for example, $\pi(0)$ is the Dynkin pyramid itself. Then,
$\Gamma(p)$ is the grading induced by declaring that each matrix
unit $e_{i,j}$ is of degree $\col(i)-\col(j)$, the notation
$\col(i)$ now denoting the column number of the $i$th box in
$\pi(p)$. Rearranging the numbers in the boxes of $\pi(p)$ so that
$\col(1)\geq \col(2)\geq\cdots\geq \col(n)$, the characteristic  of
$\Gamma(p)$ is $(c_1,\dots,c_{n-1})$ where $c_i =
\col(i)-\col(i+1).$
Finally, recall that $p \in \mathscr P_e$ defines an integral good
grading if and only if $p$ lies on as many affine hyperplanes as the
origin. For every $p \in \mathscr P_e$, there is a point $p' \in
\mathscr P_e$ lying in the closure of the alcove containing $p$,
such that $\Gamma(p')$ is an integral good grading. This means that
in type $A$, one can always restrict attention just to integral good
gradings without losing any generality. Moreover, every integral
good grading is adjacent to an even good grading, as was noted
already in the introduction of \cite{brukle}. These nice things
definitely do not usually happen in other types.

\begin{example}\rm
Let $\lambda = (3,3,2)$ as above.
The set $\mathscr P_e$ consists of all
$p=(p_1,p_2,p_3)$ with
$3p_1+3p_2+2p_3 = 0$, $|p_1-p_2| < 1, |p_2-p_3| < 2$ and
$|p_1-p_3| < 2$:
$$
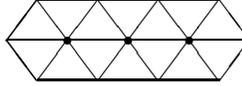
\begin{picture}(80.5,32)
\put(0,15.11){\line(1,0){92}}
\put(11.5,0){\line(1,0){69}}
\put(11.5,0){\line(1,0){69}}
\put(11.5,30.22){\line(1,0){69}}
\put(11.5,30.22){\line(1,0){69}}
\put(0,15.11){\line(3,4){11.3}}
\put(0,15.11){\line(3,4){11.3}}
\put(11.5,0){\line(3,4){22.4}}
\put(34.5,0){\line(3,4){22.4}}
\put(57.5,0){\line(3,4){22.4}}
\put(80.5,0){\line(3,4){11.2}}
\put(0,15.11){\line(3,-4){11.3}}
\put(0,15.11){\line(3,-4){11.3}}
\put(11.5,30.22){\line(3,-4){22.4}}
\put(34.5,30.22){\line(3,-4){22.4}}
\put(57.5,30.22){\line(3,-4){22.4}}
\put(80.5,30.22){\line(3,-4){11.3}}
\put(23,15.11){\circle*{3}}
\put(46,15.11){\circle*{3}}
\put(69,15.11){\circle*{3}}
\end{picture}
$$
The Weyl group $W_e \cong S_2$ is generated by the reflection in the
horizontal axis in this picture. The alcoves in $\mathscr P_e$ are
the interiors of the 14 triangles. There are just three points that
lie on as many affine hyperplanes as the origin, with associated
pyramids (renumbered so that we can read off their characteristics):
$$
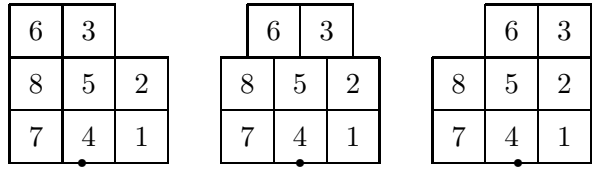
\begin{picture}(220,60)
\put(0,0){\line(1,0){60}} \put(0,20){\line(1,0){60}}
\put(0,40){\line(1,0){60}} \put(0,60){\line(1,0){40}}
\put(60,0){\line(0,1){40}} \put(0,0){\line(0,1){40}}
\put(20,0){\line(0,1){40}} \put(40,0){\line(0,1){40}}
\put(0,40){\line(0,1){20}} \put(20,40){\line(0,1){20}}
\put(40,40){\line(0,1){20}} \put(27.5,0){\circle*{3}}
\put(10,10){\makebox(0,0){{7}}} \put(30,10){\makebox(0,0){{4}}}
\put(50,10){\makebox(0,0){{1}}} \put(10,30){\makebox(0,0){{8}}}
\put(30,30){\makebox(0,0){{5}}} \put(50,30){\makebox(0,0){{2}}}
\put(10,50){\makebox(0,0){{6}}} \put(30,50){\makebox(0,0){{3}}}

\put(80,0){\line(1,0){60}} \put(80,20){\line(1,0){60}}
\put(80,40){\line(1,0){60}} \put(90,60){\line(1,0){40}}
\put(140,0){\line(0,1){40}} \put(80,0){\line(0,1){40}}
\put(100,0){\line(0,1){40}} \put(120,0){\line(0,1){40}}
\put(90,40){\line(0,1){20}} \put(110,40){\line(0,1){20}}
\put(130,40){\line(0,1){20}} \put(110,0){\circle*{3}}
\put(90,10){\makebox(0,0){{7}}} \put(110,10){\makebox(0,0){{4}}}
\put(130,10){\makebox(0,0){{1}}} \put(90,30){\makebox(0,0){{8}}}
\put(110,30){\makebox(0,0){{5}}} \put(130,30){\makebox(0,0){{2}}}
\put(100,50){\makebox(0,0){{6}}} \put(120,50){\makebox(0,0){{3}}}

\put(160,0){\line(1,0){60}} \put(160,20){\line(1,0){60}}
\put(160,40){\line(1,0){60}} \put(180,60){\line(1,0){40}}
\put(220,0){\line(0,1){40}} \put(160,0){\line(0,1){40}}
\put(180,0){\line(0,1){40}} \put(200,0){\line(0,1){40}}
\put(180,40){\line(0,1){20}} \put(200,40){\line(0,1){20}}
\put(220,40){\line(0,1){20}} \put(192.5,0){\circle*{3}}
\put(170,10){\makebox(0,0){{7}}} \put(190,10){\makebox(0,0){{4}}}
\put(210,10){\makebox(0,0){{1}}} \put(170,30){\makebox(0,0){{8}}}
\put(190,30){\makebox(0,0){{5}}} \put(210,30){\makebox(0,0){{2}}}
\put(190,50){\makebox(0,0){{6}}} \put(210,50){\makebox(0,0){{3}}}
\end{picture}
$$
Hence, there are three conjugacy classes of
integral good gradings for $e$, with
characteristics $(0,2,0,0,2,0,0)$, $(0,1,1,0,1,1,0)$ and $(0,0,2,0,0,2,0)$,
respectively.
\end{example}

\section{Good gradings for $\sp_{2n}(\C)$} \label{sp_2n}

Next, we discuss
$\g = \sp_{2n}(\C)$.
Let $V$ denote the natural $2n$-dimensional $\g$-module
with standard basis $v_1,\dots,v_n, v_{-n},\dots,v_{-1}$
and $\g$-invariant
skew-symmetric bilinear form $(.,.)$ defined by
$(v_i, v_j) = (v_{-i},v_{-j}) = 0$ and $(v_i, v_{-j}) = \delta_{i,j}$
for $1 \leq i,j \leq n$.
Let $\t$ be the set of all elements of $\g$ which act diagonally
on the standard basis of $V$.
For the simple roots $\Delta \subset \Phi$, we take
$\alpha_1 = \delta_1-\delta_2,\dots,\alpha_{n-1}=\delta_{n-1}-\delta_n$,
$\alpha_n = 2 \delta_n$, where $\delta_i \in \mathfrak{h}^*$
is defined from
$h v_i = \delta_i(h) v_i$ for each $h \in \mathfrak{h}$
and $i=1,\dots,n$.
Writing $e_{i,j}$ for the $ij$-matrix unit,
the following matrices give a Chevalley basis for $\g$:
$$
\{e_{i,j} - e_{-j,-i}\}_{1 \leq i,j \leq n}
\cup
\{e_{i,-j}+e_{j,-i}, e_{-i,j}+e_{-j,i}\}_{1 \leq i<j \leq n}
\cup
\{e_{k,-k}, e_{-k,k}\}_{1 \leq k \leq n}.
$$
Let $\sigma_{i,j} \in \{\pm 1\}$
denote the $e_{i,j}$-coefficient of the unique element in the above
basis that involves $e_{i,j}$.

Nilpotent orbits in $\g$ are parameterized by partitions $\lambda =
(\lambda_1\geq\lambda_2 \geq \cdots)$ of $2n$ such that every odd
part appears with even multiplicity. Fix such a symplectic partition
$\lambda$ throughout the section. We begin by introducing the {\em
symplectic Dynkin pyramid} of shape $\lambda$, following the idea of
\cite{elakac} closely once more. This is a diagram consisting of
$2n$ boxes each of size $2$ units by $2$ units drawn in the
$xy$-plane. As before, the coordinates of a box are the coordinates
of its midpoint, and the row and column numbers of a box mean its
$y$- and $x$-coordinate, respectively. Before we attempt a formal
definition, here are some examples of symplectic Dynkin pyramids for
$\lambda = (2,2,1,1),(4,3,3,2,2),(4,2,1,1)$ and $(4,2,2,2)$,
respectively:
$$
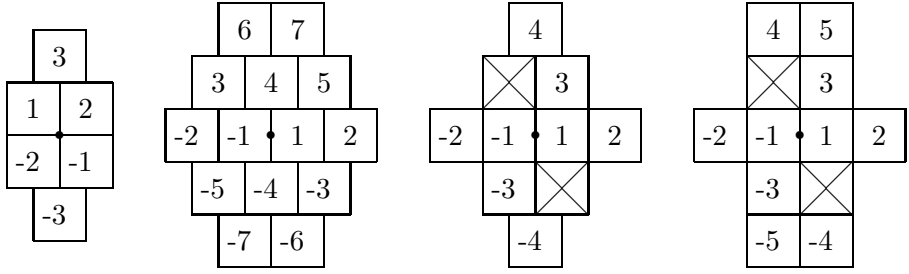
\begin{picture}(340,100)
\put(10,10){\line(1,0){20}}
\put(0,30){\line(1,0){40}}
\put(0,50){\line(1,0){40}}
\put(0,70){\line(1,0){40}}
\put(10,90){\line(1,0){20}}
\put(10,10){\line(0,1){20}}
\put(30,10){\line(0,1){20}}
\put(10,90){\line(0,-1){20}}
\put(30,90){\line(0,-1){20}}
\put(0,30){\line(0,1){40}}
\put(20,30){\line(0,1){40}}
\put(40,30){\line(0,1){40}}
\put(20,50){\circle*{3}}
\put(20,80){\makebox(0,0){{3}}}
\put(18,20){\makebox(0,0){{-3}}}
\put(8,40){\makebox(0,0){{-2}}}
\put(28,40){\makebox(0,0){{-1}}}
\put(10,60){\makebox(0,0){{1}}}
\put(30,60){\makebox(0,0){{2}}}

\put(80,0){\line(1,0){40}} \put(70,20){\line(1,0){60}}
\put(60,40){\line(1,0){80}} \put(60,60){\line(1,0){80}}
\put(70,80){\line(1,0){60}} \put(80,100){\line(1,0){40}}
\put(80,0){\line(0,1){20}} \put(100,0){\line(0,1){20}}
\put(120,0){\line(0,1){20}} \put(80,100){\line(0,-1){20}}
\put(100,100){\line(0,-1){20}} \put(120,100){\line(0,-1){20}}
\put(60,40){\line(0,1){20}} \put(80,40){\line(0,1){20}}
\put(100,40){\line(0,1){20}} \put(120,40){\line(0,1){20}}
\put(140,40){\line(0,1){20}} \put(70,20){\line(0,1){20}}
\put(90,20){\line(0,1){20}} \put(110,20){\line(0,1){20}}
\put(130,20){\line(0,1){20}} \put(70,80){\line(0,-1){20}}
\put(90,80){\line(0,-1){20}} \put(110,80){\line(0,-1){20}}
\put(130,80){\line(0,-1){20}} \put(100,50){\circle*{3}}
\put(88,10){\makebox(0,0){{-7}}} \put(108,10){\makebox(0,0){{-6}}}
\put(78,30){\makebox(0,0){{-5}}} \put(98,30){\makebox(0,0){{-4}}}
\put(118,30){\makebox(0,0){{-3}}} \put(90,90){\makebox(0,0){{6}}}
\put(110,90){\makebox(0,0){{7}}} \put(80,70){\makebox(0,0){{3}}}
\put(100,70){\makebox(0,0){{4}}} \put(120,70){\makebox(0,0){{5}}}
\put(68,50){\makebox(0,0){{-2}}} \put(88,50){\makebox(0,0){{-1}}}
\put(110,50){\makebox(0,0){{1}}} \put(130,50){\makebox(0,0){{2}}}
\put(190,0){\line(1,0){20}}
\put(190,0){\line(0,1){20}}
\put(210,0){\line(0,1){20}}
\put(190,80){\line(0,1){20}}
\put(210,80){\line(0,1){20}}
\put(190,100){\line(1,0){20}}
\put(180,20){\line(1,0){40}}
\put(180,20){\line(0,1){20}}
\put(200,20){\line(0,1){20}}
\put(220,20){\line(0,1){20}}
\put(180,80){\line(1,0){40}}
\put(180,60){\line(0,1){20}}
\put(200,60){\line(0,1){20}}
\put(220,60){\line(0,1){20}}
\put(160,40){\line(1,0){80}}
\put(160,60){\line(1,0){80}}
\put(160,40){\line(0,1){20}}
\put(180,40){\line(0,1){20}}
\put(200,40){\line(0,1){20}}
\put(220,40){\line(0,1){20}}
\put(240,40){\line(0,1){20}}
\put(200,50){\circle*{3}}
\put(198,10){\makebox(0,0){{-4}}}
\put(188,30){\makebox(0,0){{-3}}}
\put(168,50){\makebox(0,0){{-2}}}
\put(188,50){\makebox(0,0){{-1}}}
\put(210,50){\makebox(0,0){{1}}}
\put(230,50){\makebox(0,0){{2}}}
\put(210,70){\makebox(0,0){{3}}}
\put(180,60){\line(1,1){20}}
\put(180,80){\line(1,-1){20}}
\put(200,20){\line(1,1){20}}
\put(200,40){\line(1,-1){20}}
\put(200,90){\makebox(0,0){{4}}}
\put(280,0){\line(1,0){40}}
\put(280,0){\line(0,1){20}}
\put(300,0){\line(0,1){20}}
\put(320,0){\line(0,1){20}}
\put(280,60){\line(0,1){20}}
\put(280,80){\line(0,1){20}}
\put(300,80){\line(0,1){20}}
\put(320,80){\line(0,1){20}}
\put(280,100){\line(1,0){40}}
\put(280,20){\line(1,0){40}}
\put(280,20){\line(0,1){20}}
\put(300,20){\line(0,1){20}}
\put(320,20){\line(0,1){20}}
\put(280,80){\line(1,0){40}}
\put(300,60){\line(0,1){20}}
\put(320,60){\line(0,1){20}}
\put(260,40){\line(1,0){80}}
\put(260,60){\line(1,0){80}}
\put(260,40){\line(0,1){20}}
\put(280,40){\line(0,1){20}}
\put(300,40){\line(0,1){20}}
\put(320,40){\line(0,1){20}}
\put(340,40){\line(0,1){20}}
\put(300,50){\circle*{3}}
\put(288,10){\makebox(0,0){{-5}}}
\put(288,30){\makebox(0,0){{-3}}}
\put(308,10){\makebox(0,0){{-4}}}
\put(268,50){\makebox(0,0){{-2}}}
\put(288,50){\makebox(0,0){{-1}}}
\put(310,50){\makebox(0,0){{1}}}
\put(330,50){\makebox(0,0){{2}}}
\put(310,70){\makebox(0,0){{3}}}
\put(290,90){\makebox(0,0){{4}}}
\put(310,90){\makebox(0,0){{5}}}
\put(280,60){\line(1,1){20}}
\put(280,80){\line(1,-1){20}}
\put(300,20){\line(1,1){20}}
\put(300,40){\line(1,-1){20}}
\end{picture}
$$
To describe the Dynkin pyramid in the general case, the parts of
$\lambda$ indicate the number of boxes in each row, and the rows are
added to the diagram in order, starting with the row corresponding
to the largest part of $\lambda$ closest to the $x$-axis and moving
out from there, in a centrally symmetric way. The only complication
is that if some (necessarily even) part $\lambda_i$ of $\lambda$ has
odd multiplicity, then the first time a row of this length is added
to the diagram it is split into two halves, the right half is added
to the next free row in the upper half plane in columns
$1,3,\dots,\lambda_i-1$ and the left half is added to the lower half
plane in a centrally symmetric way. We refer to the exceptional rows
arising in this way as {\em skew rows}; in particular, if the
largest part of $\lambda$ has odd multiplicity, then the zeroth row
is a skew row.  The missing boxes in skew rows are drawn as box with
a cross through it.  We let $r_1 < \cdots < r_m$ denote the numbers
of the non-empty rows in the upper half plane that are not skew
rows, and define $\bar \lambda_i$ to be the number of boxes in row
$r_i$ for each $i=1,\dots,m$.

Fix from now on a numbering of the boxes of the Dynkin pyramid by
the numbers $1,\dots,n,-n,\dots,-1$ in such a way that $i$ and $-i$
appear in centrally symmetric boxes, for each $i=1,\dots,n$.
\iffalse
; and so
that the positive labels occur in the positive half-plane ($x \ge 0$
and $y \ge 0$ if $x = 0$).
\fi
As before, we write $\row(i)$ and
$\col(i)$ for the row and column numbers of the $i$th box. Now we
can fix a choice of an $\sl_2$-triple $(e,h,f)$ with $e$ of Jordan
type $\lambda$. Define $e \in \g$ to be the matrix $\sum_{i,j}
\sigma_{i,j} e_{i,j}$, where the sum is over all pairs $i,j$ of
boxes in the Dynkin pyramid such that
\begin{itemize}
\item[{\em either}\:] $\row(i) = \row(j)$
and $\col(i) = \col(j)+2$;
\item[{\em or}\:] $\row(i) = -\row(j)$
is a skew row in the upper half plane
and $\col(i) = 1, \col(j) = -1$.
\end{itemize}
For example if $\lambda = (4,2,1,1)$ and the Dynkin pyramid
is labelled as above, then
$e = e_{3,-3}+e_{2,1}+e_{1,-1}-e_{-1,-2}$.
Also define $h \in \mathfrak{h}$ to be
$\sum_i \col(i) e_{i,i}$, again summing over all boxes in the Dynkin
pyramid. There is then a unique $f \in \g$
such that $(e,h,f)$ is an $\sl_2$-triple.

The important thing about these choices is that once again
 $\t$ is contained in $\c$ and $\t_e$
is a Cartan subalgebra of $\c_e$. In fact, $E_e$ consists of all
matrices in $\t$ with entries from $\R$, such that the $i$th and
$j$th diagonal entries are equal whenever $\row(i) = \row(j)$ and
the $k$th diagonal entry is zero whenever $\row(k)$ is a skew row,
for $1 \leq i,j,k \leq n$. For $i=1,\dots,m$, let $\eps_i$ be the
function picking out the $j$th diagonal entry of a matrix in $E_e$,
where $1 \leq j \leq n$ here is chosen so that the $j$th box is in
row $r_i$ in the Dynkin pyramid. Then, $\eps_1,\dots,\eps_m$ form a
basis for $E_e^*$. We identify $E_e$ and $E_e^*$ via the trace form
$(\eps_i,\eps_j) = \frac{\delta_{i,j}}{2\bar \lambda_i}$. We have
that
\begin{align*}
\Phi_e &= \{\eps_i \pm \eps_j, \pm 2\eps_k \mid 1 \leq i,j,k \leq m,
i \neq j\}\\\intertext{if there are no skew rows, i.e.\ all non-zero
parts of $\lambda$ are of even multiplicity, or} \Phi_e &= \{\pm
\eps_h, \eps_i \pm \eps_j, \pm 2\eps_k \mid 1 \leq h,i, j, k \leq m,
i \neq j\}
\end{align*}
if there are skew rows.
The values of $d(\alpha)$ from Theorem~\ref{ggp} 
for all $\alpha \in \Phi_e$ are:
\begin{align*}
d(\eps_i \pm \eps_j) & = 1 + |\bar\lambda_i - \bar\lambda_j|, \\
d(\pm 2\eps_k) & = \left\{ \begin{array}{ll} 1 & \text{if $\bar\lambda_k$ is odd} \\
                                 3 & \text{if $\bar\lambda_k$ is
                                 even,} \end{array} \right. \\
d(\pm\eps_h) & = 1 + \min\{|\bar \lambda_h - t| \mid
\text{$t$ is a non-zero part of $\lambda$ of odd multiplicity}\}.
\end{align*}
Hence, representing a point $p \in E_e$ as an
$m$-tuple $p=(p_1,\dots,p_m)$ of real numbers defined from $p_i =
\eps_i(p)$,
the good
grading polytope $\mathscr P_e$ is the open subset of $E_e$ defined
by the inequalities
\begin{align*}
|p_i \pm p_j| &< 1 + \bar\lambda_i - \bar\lambda_j,\\
|p_k|&< \left\{
\begin{array}{ll}
\frac{1}{2}&\text{if $\bar \lambda_k$ is odd,}\\
1&\text{if $\bar \lambda_k$ is even of multiplicity $> 2$ in $\lambda$,}\\
\frac{3}{2}&\text{if $\bar \lambda_k$ is even of multiplicity 2 in $\lambda$,}
\end{array}\right.
\end{align*}
for all $1 \leq i < j \leq m$ and $1 \leq k \leq m$.
Also, letting $\overline m_i$ denote the multiplicity of $i$ as a part
of the partition $(\bar\lambda_1,\bar\lambda_2,\dots)$,
the restricted Weyl group $W_e$
is the subgroup
$B_{\overline m_1} \times B_{\overline m_2} \times \cdots$
of the Weyl group $B_m$,
acting on $\{\pm\eps_1,\dots,\pm \eps_m\}$ by all sign changes and all
permutations
of $\eps_i$'s of equal length.

Again, there is a useful combinatorial way to visualize the good
grading $\Gamma(p)$ corresponding to a point $p=(p_1,\dots,p_m) \in
\mathscr P_e$ involving pyramids. Let $\pi(p)$ denote the pyramid
obtained from the symplectic Dynkin pyramid by sliding all numbered
boxes in rows $\pm r_i$ to the right by $\pm p_i$ units. Then,
$\Gamma(p)$ is the grading induced by declaring that each matrix
unit $e_{i,j}$ is of degree $\col(i)-\col(j)$, where the notation
$\col(i)$ now denotes the column number of the $i$th box in
$\pi(p)$. The characteristic of the grading $\Gamma(p)$ can be
computed by first rearranging the entries in the boxes of $\pi(p)$
using all permutations and sign changes from the Weyl group $W=B_n$,
so that $\col(1)\geq \col(2)\geq\cdots\geq \col(n)\geq 0$. Then, the
characteristic of $\Gamma(p)$ is $(c_1,\dots,c_n)$ where $c_i =
\col(i) - \col(i+1)$ for $i=1,\dots,n-1$ and $c_n = 2 \col(n)$.

\iffalse The integral good gradings for $e$ that are compatible with
$\t$ are the gradings $\Gamma(p)$ for all $p=(p_1,\dots,p_m) \in
\mathscr P_e$ that  satisfy $p_i \pm p_j \in \Z$ and either $p_k \in
\Z$ if there are skew rows or $2p_k \in \Z$ if there are not, for
all $1 \leq i < j \leq m$ and $1 \leq k \leq m$. For every $p \in
\mathscr P_e$, the closure of the alcove containing $p$ always
contains a point $q$ such that $\Gamma(q)$ is a good
$\frac{1}{2}\Z$-grading for $e$. Every good grading for $e$ is
adjacent to an integral good grading. \iffalse Finally, there exists
an even good grading for $e$ if and only if all parts of $\lambda$
are even, or $\lambda$ has odd parts and all even parts of $\lambda$
have multiplicity 2. \fi This should be compared with the situation
in type $A$: it is no longer sufficient to restrict attention just
to integral good gradings. \fi

\begin{example}\rm
Take $\lambda = (2,2,1,1)$.
Then $\mathscr P_e$ consists of all $p = (p_1,p_2)$ with
$|p_1| < 3/2$ and $|p_2| < 1/2$:
$$
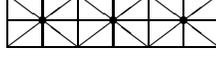
\begin{picture}(80,20)
\put(0,0){\line(1,0){80}}
\put(0,0){\line(4,3){26.66}}
\put(26.66,0){\line(4,3){26.66}}
\put(53.33,0){\line(4,3){26.66}}
\put(0,20){\line(4,-3){26.66}}
\put(26.66,20){\line(4,-3){26.66}}
\put(53.33,20){\line(4,-3){26.66}}
\put(0,10){\line(1,0){80}}
\put(0,20){\line(1,0){80}}
\put(0,0){\line(0,1){20}}
\put(13.33,0){\line(0,1){20}}
\put(26.66,0){\line(0,1){20}}
\put(40,0){\line(0,1){20}}
\put(53.33,0){\line(0,1){20}}
\put(66.66,0){\line(0,1){20}}
\put(80,0){\line(0,1){20}}
\put(13.33,10){\circle*{3}}
\put(40,10){\circle*{3}}
\put(66.66,10){\circle*{3}}
\end{picture}
$$
The Weyl group $W_e \cong S_2 \times S_2$ is generated by
reflections in the horizontal and vertical axes. There are three
integral good gradings for $e$ compatible with $\t$, with associated
pyramids (renumbered so that we can read off their characteristics):
$$
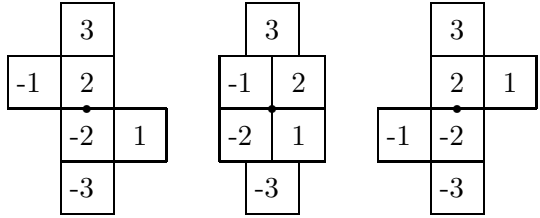
\begin{picture}(180,80)
\put(150,0){\line(1,0){20}} \put(130,20){\line(1,0){40}}
\put(130,40){\line(1,0){60}} \put(150,60){\line(1,0){40}}
\put(150,80){\line(1,0){20}} \put(150,0){\line(0,1){20}}
\put(170,0){\line(0,1){20}} \put(150,80){\line(0,-1){20}}
\put(170,80){\line(0,-1){20}} \put(130,20){\line(0,1){20}}
\put(150,20){\line(0,1){40}} \put(170,20){\line(0,1){40}}
\put(190,40){\line(0,1){20}} \put(160,40){\circle*{3}}
\put(160,70){\makebox(0,0){{3}}} \put(158,10){\makebox(0,0){{-3}}}
\put(158,30){\makebox(0,0){{-2}}} \put(180,50){\makebox(0,0){{1}}}
\put(138,30){\makebox(0,0){{-1}}} \put(160,50){\makebox(0,0){{2}}}
\put(80,0){\line(1,0){20}} \put(70,20){\line(1,0){40}}
\put(70,40){\line(1,0){40}} \put(70,60){\line(1,0){40}}
\put(80,80){\line(1,0){20}} \put(80,0){\line(0,1){20}}
\put(100,0){\line(0,1){20}} \put(80,80){\line(0,-1){20}}
\put(100,80){\line(0,-1){20}} \put(70,20){\line(0,1){40}}
\put(90,20){\line(0,1){40}} \put(110,20){\line(0,1){40}}
\put(90,40){\circle*{3}} \put(90,70){\makebox(0,0){{3}}}
\put(88,10){\makebox(0,0){{-3}}} \put(78,30){\makebox(0,0){{-2}}}
\put(100,30){\makebox(0,0){{1}}} \put(78,50){\makebox(0,0){{-1}}}
\put(100,50){\makebox(0,0){{2}}}
\put(10,0){\line(1,0){20}} \put(10,20){\line(1,0){40}}
\put(-10,40){\line(1,0){60}} \put(-10,60){\line(1,0){40}}
\put(10,80){\line(1,0){20}} \put(-10,40){\line(0,1){20}}
\put(50,20){\line(0,1){20}} \put(10,0){\line(0,1){80}}
\put(30,0){\line(0,1){80}} \put(20,40){\circle*{3}}
\put(20,70){\makebox(0,0){{3}}} \put(18,10){\makebox(0,0){{-3}}}
\put(18,30){\makebox(0,0){{-2}}} \put(40,30){\makebox(0,0){{1}}}
\put(-2,50){\makebox(0,0){{-1}}} \put(20,50){\makebox(0,0){{2}}}
\end{picture}
$$
These have characteristics $(2,0,0)$, $(0,1,0)$ and $(2,0,0)$,
respectively.  The right and left good gradings are conjugate by the
element of $W_e$ corresponding to the reflection in the vertical
axes of the good grading polytope.
\end{example}

\section{Good gradings for $\so_N(\C)$} \label{so_N}

Finally let $\g = \so_{N}(\C)$ and set $n = \lfloor
\frac{N}{2}\rfloor$, assuming $N \geq 3$. Let $V$ be the natural
$N$-dimensional $\g$-module with standard basis $v_1,\dots,v_n, v_0,
v_{-n},\dots,v_{-1}$ and $\g$-invariant symmetric bilinear form
$(.,.)$ defined by $(v_0, v_i) = (v_0, v_{-i}) = 0, (v_0,v_0) = 2$,
$(v_i, v_j) = (v_{-i},v_{-j}) = 0$ and $(v_i, v_{-j}) =
\delta_{i,j}$ for $1 \leq i,j \leq n$ (omitting $v_0$ everywhere if
$N$ is even). Let $\t$ be the set of all elements of $\g$ which act
diagonally on the standard basis of $V$. Defining $\delta_i \in
\mathfrak{h}^*$ as in section \ref{sp_2n}, a choice of simple roots
$\Delta \subset \Phi$ is given by $\alpha_1 =
\delta_1-\delta_2,\dots,\alpha_{n-1}=\delta_{n-1}-\delta_n$, and
$\alpha_n = \delta_{n-1}+\delta_n$ if $N$ is even or $\delta_n$ if
$N$ is odd. The following matrices give a Chevalley basis for $\g$
(again omitting the last family if $N$ is even):
\begin{align*}
\{e_{i,j} - e_{-j,-i}\}_{1 \leq i,j \leq n} &\cup
\{e_{i,-j}-e_{j,-i}, e_{-j,i}-e_{-i,j}\}_{1 \leq i < j \leq n}\\
&\cup \{2e_{k,0}-e_{0,-k}, e_{0,k}-2e_{-k,0}\}_{1 \leq k \leq n}.
\end{align*}
As before, define $\sigma_{i,j}$ to be the coefficient of $e_{i,j}$ in this basis if it appears, or zero if no basis element involves $e_{i,j}$.

Nilpotent orbits in $\g$ are parameterized by partitions $\lambda =
(\lambda_1\geq\lambda_2 \geq \cdots)$ of $N$ such that every even
part appears with even multiplicity; in case $N$ is even, we mean
nilpotent orbits under the group $\mathrm{O}_N$ not $\mathrm{SO}_N$
here. Fix such an orthogonal partition $\lambda$ throughout the
section. We need the {\em orthogonal Dynkin pyramid} of type
$\lambda$, which again consists of $N$ boxes of size $2$ units by
$2$ units arranged in the $xy$-plane in a centrally symmetric way.
Assume to start with that $N$ is even. Then the Dynkin pyramid is
constructed like in the symplectic case, adding rows of lengths
determined by the parts of $\lambda$ working outwards from the
$x$-axis starting with the largest part, in a centrally symmetric
way. The only difficulty is if some (necessarily odd) part of
$\lambda$ appears with odd multiplicity. As $N$ is even, the number
of distinct parts having odd multiplicity is even. Choose $i_1 < j_1
< \cdots < i_r < j_r$ such that $\lambda_{i_1} > \lambda_{j_1}
>\dots>\lambda_{i_r}>\lambda_{j_r}$ are representatives for all the
distinct odd parts of $\lambda$ having odd multiplicity. Then the
first time the part $\lambda_{i_s}$ needs to be added to the
diagram, the part $\lambda_{j_s}$ is also added at the same time, so
that the parts $\lambda_{i_s}$ and $\lambda_{j_s}$ of $\lambda$
contribute two centrally symmetric rows to the diagram, one row in
the upper half plane with boxes in columns
$1-\lambda_{j_s},3-\lambda_{j_s},\dots,\lambda_{i_s}-1$ and the
other row in the lower half plane with boxes in columns
$1-\lambda_{i_s},3-\lambda_{i_s},\dots,\lambda_{j_s}-1$. We will
refer to the exceptional rows arising in this way as {\em skew
rows}. Here are some examples, for $\lambda =
(3,3,2,2),(3,1,1,1),(3,2,2,1)$ and $(7,7,7,3)$, respectively:
$$
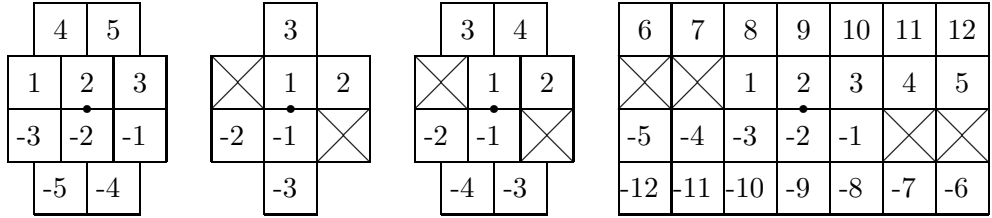
\begin{picture}(371,80)
\put(10,0){\line(1,0){40}}
\put(0,20){\line(1,0){60}}
\put(0,40){\line(1,0){60}}
\put(0,60){\line(1,0){60}}
\put(10,80){\line(1,0){40}}
\put(10,0){\line(0,1){20}}
\put(30,0){\line(0,1){20}}
\put(50,0){\line(0,1){20}}
\put(10,60){\line(0,1){20}}
\put(30,60){\line(0,1){20}}
\put(50,60){\line(0,1){20}}
\put(0,20){\line(0,1){40}}
\put(20,20){\line(0,1){40}}
\put(40,20){\line(0,1){40}}
\put(60,20){\line(0,1){40}}
\put(30,40){\circle*{3}}
\put(18,10){\makebox(0,0){{-5}}}
\put(20,70){\makebox(0,0){{4}}}
\put(8,30){\makebox(0,0){{-3}}}
\put(10,50){\makebox(0,0){{1}}}
\put(28,30){\makebox(0,0){{-2}}}
\put(30,50){\makebox(0,0){{2}}}
\put(50,50){\makebox(0,0){{3}}}
\put(48,30){\makebox(0,0){{-1}}}
\put(38,10){\makebox(0,0){{-4}}}
\put(40,70){\makebox(0,0){{5}}}
\put(97,0){\line(1,0){20}}
\put(77,20){\line(1,0){60}}
\put(77,40){\line(1,0){60}}
\put(77,60){\line(1,0){60}}
\put(97,80){\line(1,0){20}}
\put(97,0){\line(0,1){80}}
\put(77,20){\line(0,1){40}}
\put(117,0){\line(0,1){80}}
\put(137,20){\line(0,1){40}}
\put(117,20){\line(1,1){20}}
\put(77,40){\line(1,1){20}}
\put(117,40){\line(1,-1){20}}
\put(77,60){\line(1,-1){20}}
\put(107,40){\circle*{3}}
\put(105,10){\makebox(0,0){{-3}}}
\put(105,30){\makebox(0,0){{-1}}}
\put(85,30){\makebox(0,0){{-2}}}
\put(107,70){\makebox(0,0){{3}}}
\put(107,50){\makebox(0,0){{1}}}
\put(127,50){\makebox(0,0){{2}}}
\put(164,0){\line(1,0){40}}
\put(154,20){\line(1,0){60}}
\put(154,40){\line(1,0){60}}
\put(154,40){\line(1,1){20}}
\put(154,60){\line(1,-1){20}}
\put(194,20){\line(1,1){20}}
\put(194,40){\line(1,-1){20}}
\put(154,60){\line(1,0){60}}
\put(164,80){\line(1,0){40}}
\put(164,0){\line(0,1){20}}
\put(184,0){\line(0,1){20}}
\put(204,0){\line(0,1){20}}
\put(154,20){\line(0,1){40}}
\put(174,20){\line(0,1){40}}
\put(194,20){\line(0,1){40}}
\put(214,20){\line(0,1){40}}
\put(164,60){\line(0,1){20}}
\put(184,60){\line(0,1){20}}
\put(204,60){\line(0,1){20}}
\put(184,40){\circle*{3}}
\put(192,10){\makebox(0,0){{-3}}}
\put(194,70){\makebox(0,0){{4}}}
\put(204,50){\makebox(0,0){{2}}}
\put(162,30){\makebox(0,0){{-2}}}
\put(184,50){\makebox(0,0){{1}}}
\put(182,30){\makebox(0,0){{-1}}}
\put(172,10){\makebox(0,0){{-4}}}
\put(174,70){\makebox(0,0){{3}}}
\put(231,0){\line(1,0){140}}
\put(231,20){\line(1,0){140}}
\put(231,40){\line(1,0){140}}
\put(231,60){\line(1,0){140}}
\put(231,80){\line(1,0){140}}
\put(231,00){\line(0,1){80}}
\put(251,0){\line(0,1){80}}
\put(271,0){\line(0,1){80}}
\put(291,0){\line(0,1){80}}
\put(231,40){\line(1,1){20}}
\put(351,40){\line(1,-1){20}}
\put(351,20){\line(1,1){20}}
\put(231,60){\line(1,-1){20}}
\put(251,40){\line(1,1){20}}
\put(331,40){\line(1,-1){20}}
\put(331,20){\line(1,1){20}}
\put(251,60){\line(1,-1){20}}
\put(311,0){\line(0,1){80}}
\put(331,0){\line(0,1){80}}
\put(351,0){\line(0,1){80}}
\put(371,0){\line(0,1){80}}
\put(301,40){\circle*{3}}
\put(239,10){\makebox(0,0){{-12}}}
\put(241,70){\makebox(0,0){{6}}}
\put(281,50){\makebox(0,0){{1}}}
\put(281,70){\makebox(0,0){{8}}}
\put(301,70){\makebox(0,0){{9}}}
\put(321,70){\makebox(0,0){{10}}}
\put(321,50){\makebox(0,0){{3}}}
\put(341,50){\makebox(0,0){{4}}}
\put(341,70){\makebox(0,0){{11}}}
\put(361,70){\makebox(0,0){{12}}}
\put(261,70){\makebox(0,0){{7}}}
\put(361,50){\makebox(0,0){{5}}}
\put(259,10){\makebox(0,0){{-11}}}
\put(279,10){\makebox(0,0){{-10}}}
\put(239,30){\makebox(0,0){{-5}}}
\put(259,30){\makebox(0,0){{-4}}}
\put(279,30){\makebox(0,0){{-3}}}
\put(299,30){\makebox(0,0){{-2}}}
\put(299,10){\makebox(0,0){{-9}}}
\put(301,50){\makebox(0,0){{2}}}
\put(339,10){\makebox(0,0){{-7}}}
\put(319,30){\makebox(0,0){{-1}}}
\put(319,10){\makebox(0,0){{-8}}}
\put(359,10){\makebox(0,0){{-6}}}
\end{picture}
$$
If $N$ is odd, there is one additional consideration.
There must be some odd part appearing with odd multiplicity.
Let $\lambda_i$ be the largest such part, and put
$\lambda_i$ boxes into the zeroth row
in columns $1-\lambda_i,3-\lambda_i,\dots,\lambda_i-1$;
we also treat this zeroth row as a {skew row}.
Now remove the part $\lambda_i$
from $\lambda$,
to obtain a partition of an even number. The remaining parts are then added
to the diagram exactly as in the case $N$ even.
We give two more examples, for $\lambda = (6,6,5)$ and $(5,3,1)$, respectively:
$$
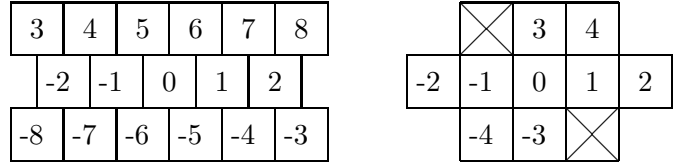
\begin{picture}(250,60)
\put(0,0){\line(1,0){120}}
\put(0,20){\line(1,0){120}}
\put(0,40){\line(1,0){120}}
\put(0,60){\line(1,0){120}}
\put(0,0){\line(0,1){20}}
\put(20,0){\line(0,1){20}}
\put(40,0){\line(0,1){20}}
\put(60,0){\line(0,1){20}}
\put(80,0){\line(0,1){20}}
\put(100,0){\line(0,1){20}}
\put(120,0){\line(0,1){20}}
\put(0,40){\line(0,1){20}}
\put(20,40){\line(0,1){20}}
\put(40,40){\line(0,1){20}}
\put(60,40){\line(0,1){20}}
\put(80,40){\line(0,1){20}}
\put(100,40){\line(0,1){20}}
\put(120,40){\line(0,1){20}}
\put(10,20){\line(0,1){20}}
\put(30,20){\line(0,1){20}}
\put(50,20){\line(0,1){20}}
\put(70,20){\line(0,1){20}}
\put(90,20){\line(0,1){20}}
\put(110,20){\line(0,1){20}}
\put(60,30){\makebox(0,0){{0}}}
\put(8,10){\makebox(0,0){{-8}}}
\put(28,10){\makebox(0,0){{-7}}}
\put(48,10){\makebox(0,0){{-6}}}
\put(68,10){\makebox(0,0){{-5}}}
\put(88,10){\makebox(0,0){{-4}}}
\put(108,10){\makebox(0,0){{-3}}}
\put(18,30){\makebox(0,0){{-2}}}
\put(38,30){\makebox(0,0){{-1}}}
\put(10,50){\makebox(0,0){{3}}}
\put(80,30){\makebox(0,0){{1}}}
\put(100,30){\makebox(0,0){{2}}}
\put(30,50){\makebox(0,0){{4}}}
\put(50,50){\makebox(0,0){{5}}}
\put(70,50){\makebox(0,0){{6}}}
\put(90,50){\makebox(0,0){{7}}}
\put(110,50){\makebox(0,0){{8}}}
\put(170,0){\line(1,0){60}}
\put(150,20){\line(1,0){100}}
\put(150,40){\line(1,0){100}}
\put(170,60){\line(1,0){60}}
\put(170,60){\line(1,-1){20}}
\put(210,0){\line(1,1){20}}
\put(210,20){\line(1,-1){20}}
\put(170,40){\line(1,1){20}}
\put(150,20){\line(0,1){20}}
\put(170,0){\line(0,1){60}}
\put(190,0){\line(0,1){60}}
\put(210,0){\line(0,1){60}}
\put(230,0){\line(0,1){60}}
\put(250,20){\line(0,1){20}}
\put(200,30){\makebox(0,0){{0}}}
\put(178,10){\makebox(0,0){{-4}}}
\put(198,10){\makebox(0,0){{-3}}}
\put(178,30){\makebox(0,0){{-1}}}
\put(200,50){\makebox(0,0){{3}}}
\put(220,50){\makebox(0,0){{4}}}
\put(240,30){\makebox(0,0){{2}}}
\put(158,30){\makebox(0,0){{-2}}}
\put(220,30){\makebox(0,0){{1}}}
\end{picture}
$$
Let $r_1 < \cdots < r_m$ denote the numbers of the non-empty rows
in the upper half plane of the Dynkin pyramid
that are not skew rows, and define
$\bar \lambda_i$ to be the number of boxes in row $r_i$ for each
$i=1,\dots,m$.

Note in the case $N$ is odd,
there is always a box at $(0,0)$; we always number it by $0$. The remaining boxes, for $N$ even or odd, are numbered $\pm 1.\dots,\pm n$
exactly as in the symplectic case, and we use the
notation $\row(i)$ and $\col(i)$ just as before.
Define $e \in \g$ to be the matrix $\sum_{i,j} \sigma_{i,j}
e_{i,j}$, where the sum is over all pairs $i,j$ of boxes in the
Dynkin pyramid such that
\begin{itemize}
\item[{\em either}\:] $\row(i) = \row(j)$
and $\col(i) = \col(j)+2$;
\item[{\em or}\:] $\row(i) = -\row(j)$
is a skew row in the upper half plane
and $\col(i) = 2, \col(j) = 0$;
\item[{\em or}\:] $\row(i) = -\row(j)$
is a skew row in the upper half plane
and $\col(i) = 0, \col(j) = -2$.
\end{itemize}
This is an element of $\g$ having Jordan type $\lambda$. (If {\em
all} parts of $\lambda$ are even then there is another conjugacy
class of elements of $\g$ of Jordan type $\lambda$, a representative
for which can be obtained using the above formula by swapping the
entries $i$ and $-i$ in the Dynkin pyramid for some $1 \leq i \leq
n$.) For example, for $\lambda = (5,3,1)$ labelled as above, $e =
e_{4,3}-e_{4,-3} +e_{3,-4} -e_{-3,-4}+ e_{2,1} + e_{1,0} - e_{0,-1}
- e_{-1,-2}$. Let $h = \sum_i \col(i) e_{i,i}$. Then there is a
unique element $f \in \g$ such that $(e,h,f)$ is an $\sl_2$-triple.

Again, these choices ensure that $\t$ is contained in $\c$ and
$\t_e$
is a Cartan subalgebra of $\c_e$.
Like for $\sp_{2n}(\C)$,
$E_e$ consists of
all matrices in $\t$ with entries from $\R$,
such that the $i$th and $j$th diagonal entries
are equal whenever $\row(i) = \row(j)$
and the $k$th
diagonal entry is zero whenever $\row(k)$ is a skew row,
for $1 \leq i,j,k \leq n$.
Define the basis $\eps_1,\dots,\eps_m$
for $E_e^*$ just as in the symplectic case,
and work with the inner product defined by
$(\eps_i,\eps_j) = \frac{\delta_{i,j}}{2\bar\lambda_i}$.
This time, we have that
\begin{align*}
\Phi_e &= \{\eps_i \pm \eps_j, \pm 2\eps_k \mid 1 \leq i,j,k \leq m,
i \neq j, \bar\lambda_k \neq 1\}\\\intertext{if there are no skew
rows, i.e.\ all non-zero parts of $\lambda$ are of even
multiplicity, or} \Phi_e &= \{\pm \eps_h, \eps_i \pm \eps_j, \pm
2\eps_k \mid 1 \leq h,i, j, k \leq m, i \neq j, \bar\lambda_k \neq
1\}
\end{align*}
if there are skew rows.
The values of $d(\alpha)$ from Theorem \ref{ggp} for all $\alpha \in \Phi_e$ 
are:
\begin{align*}
d(\eps_i \pm \eps_j) & = 1 + |\bar\lambda_i - \bar\lambda_j|, \\
d(\pm 2\eps_k) & = \left\{ \begin{array}{ll} 1 & \text{if $\bar\lambda_k$ is even} \\
                                 3 & \text{if $\bar\lambda_k$ is
                                 odd,} \end{array} \right. \\
d(\pm \eps_h) & = 1 + \min\{|\bar \lambda_h - t|
\mid \text{$t$ is a non-zero part of $\lambda$ of odd multiplicity}\}.
\end{align*}
Hence, representing points $p \in E_e$ as $m$-tuples
$p=(p_1,\dots,p_m)$ of real numbers so that $p_i = \eps_i(p)$,
the good grading polytope $\mathscr P_e$ is the open subset
of $E_e$ defined by the inequalities
\begin{align*}
|p_i \pm p_j| &< 1 + \bar\lambda_i - \bar\lambda_j,\\
|p_k|&< \left\{
\begin{array}{ll}
\frac{1}{2}&\text{if $\bar \lambda_k$ is even,}\\
1&\text{if $\bar \lambda_k$ is odd of multiplicity $> 2$ in $\lambda$,}\\
\frac{3}{2}&\text{if $\bar \lambda_k \neq 1$
is odd of multiplicity 2 in $\lambda$,}\\
s&\text{if $\bar \lambda_k = 1$
is of multiplicity 2 in $\lambda$,}
\end{array}\right.
\end{align*}
for all $1 \leq i < j \leq m$ and $1 \leq k \leq m$.
Here, in the case that the part $1$ is of multiplicity 2 in $\lambda$,
$s$ denotes the smallest part of $\lambda$ that is greater
than $1$.
Let $\overline m_i$ denote the multiplicity of $i$ as a part
of the partition $(\bar\lambda_1,\bar\lambda_2,\dots)$.
If there are skew rows, then
the restricted Weyl group $W_e$
is
$B_{\overline m_1} \times B_{\overline m_2} \times \cdots$,
acting on
$\{\pm\eps_1,\dots,\pm \eps_m\}$ by all sign changes and all permutations
of $\eps_i$'s of equal length.
If there are no skew rows, then $W_e$
is instead the subgroup of
$B_{\overline m_1} \times B_{\overline m_2} \times \cdots$
consisting of all the elements in this group that act on
$\{\pm\eps_1,\dots,\pm \eps_m\}$ with only
an even number of sign changes of the $\eps_i$'s for which $\bar\lambda_i$
is odd.

For $p=(p_1,\dots,p_m) \in \mathscr P_e$, we define the pyramid
$\pi(p)$ by sliding all numbered boxes in rows $\pm r_i$ of the
orthogonal Dynkin pyramid to the right by $\pm p_i$ units. Then, the
grading $\Gamma(p)$ associated to the point $p \in \mathscr P_e$ is
the grading induced by declaring that each matrix unit $e_{i,j}$ is
of degree $\col(i)-\col(j)$, where the notation $\col(i)$ here
denotes the column number of the $i$th box in $\pi(p)$. To compute
the characteristic of the grading $\Gamma(p)$, suppose first that
$N$ is odd. Rearrange the entries in the boxes of $\pi(p)$ using all
permutations and sign changes from the Weyl group $W=B_n$ so that
$\col(1)\geq \col(2)\geq\cdots\geq \col(n) \geq 0$. Then, the
characteristic of $\Gamma(p)$ is $(c_1,\dots,c_n)$ where $c_i =
\col(i) - \col(i+1)$ for $i=1,\dots,n-1$ and $c_n = \col(n)$.
Instead, if $N$ is even, rearrange the entries in the boxes of
$\pi(p)$ using all permutations and sign changes from the Weyl group
$W=D_n$ (i.e.\ so that there are only an even number of sign changes
in total) so that $\col(1)  \geq \cdots \geq \col(n-1) \geq
|\col(n)|$. Then, the characteristic of $\Gamma(p)$ is
$(c_1,\dots,c_n)$ where $c_i = \col(i) - \col(i+1)$ for
$i=1,\dots,n-1$ and $c_n = \col(n-1)+\col(n)$.

\ifex@
\section{Exceptional Lie algebras}

Suppose $\g$ is of exceptional type and
work with the Bala-Carter parametrization of the nilpotent
orbits. As we have said, it is routine in each case to compute
the restricted root system and the good grading polytope explicitly,
and we have implemented this in {\sc GAP} \cite{GAP4}.
However there is too much information
to record explicitly here.
Instead, we restrict our attention to integral good gradings, as classified
in terms of their characteristics in \cite{elakac}, and
list here the {\em adjacency graphs} for all integral good gradings
for all exceptional Lie algebras.
These are the graphs with vertices being the characteristics of the
integral good gradings,
with an edge joining each pair of adjacent good characteristics.
Actually, in $E_8, F_4$ and $G_2$, there are no edges at all, so there is nothing further to be said. The adjacency graphs having more than one edge
for $E_6$ and $E_7$ are listed in Tables~\ref{t1} and \ref{t2}.
One sees at once from these graphs and the tables in \cite{C,elakac}
that there are only seven nilpotent orbits for which the adjacency
graphs of integral good gradings are not connected, namely,
$\g = E_6$ and $e$ of Bala-Carter label $A_4$ or $A_4+A_1$,
$\g = E_7$ and $e$ of Bala-Carter label $A_4, A_4+A_1$ or $D_5(a_1)$,
and $\g = E_8$ and $e$ of Bala-Carter label $D_7(a_1)$ or $D_7(a_2)$.
In all of these cases, we have checked that the adjacency graphs of good $\frac{1}{2}\Z$-gradings
for $e$ are connected.
\fi

\end{document}